\documentclass[reqno,11pt]{amsart}
\usepackage{amsmath,amssymb,amsthm,a4wide,enumerate,graphicx,url,subfigure}
\usepackage[small,bf]{caption} 
\usepackage[dvipsnames]{xcolor}
\newtheorem*{remark}{Remark}
\newcommand{\prn}[1]{\left(#1\right)}
\newcommand{\abs}[1]{\left|#1\right|}

\newcommand{\pd}[2]{\frac{\partial#1}{\partial#2}}
\newcommand{\ud}[1]{\,\mathrm{d}#1}

\newcommand{\mE}{\mathcal{E}}
\newcommand{\mV}{\mathcal{V}}

\newcommand{\mR}{\mathcal{R}}
\newcommand{\LT}{\Tilde{L}}
\newcommand{\RT}{\Tilde{R}}
\newcommand{\ST}{\Tilde{S}}
\newcommand{\PT}{\Tilde{P}}
\DeclareMathOperator\erf{erf}
\newcommand{\mzero}{m_0}
\newcommand{\mone}{m_1}
\newcommand{\mtwo}{m_2}
\newcommand{\Mzero}{\mathcal{M}_0}
\newcommand{\Mone}{\mathcal{M}_1}
\newcommand{\Mtwo}{\mathcal{M}_2}
\newcommand{\dela}{\Delta a}

\title[Moment Methods for Advection on Networks and Application to Forest Pest]{Moment Methods for Advection on Networks and an Application to Forest Pest Life Cycle Models}

\author[R. Chinomona]{Rujeko Chinomona}
\address[Rujeko Chinomona]
{Department of Mathematics \\ Temple University \\ \newline
1805 North Broad Street \\ Philadelphia, PA 19122}
\email{rujeko.chinomona@temple.edu}

\author[K. Kean]{Kiera Kean}
\address[Kiera Kean]
{Department of Mathematics \\ Temple University \\ \newline
1805 North Broad Street \\ Philadelphia, PA 19122}
\email{kiera.kean@temple.edu}

\author[B. Seibold]{Benjamin Seibold}
\address[Benjamin Seibold]
{Department of Mathematics \\ Temple University \\ \newline
1805 North Broad Street \\ Philadelphia, PA 19122}
\email{seibold@temple.edu}
\urladdr{http://www.math.temple.edu/\~{}seibold}

\author[J. Woods]{Jacob Woods}
\address[Jacob Woods]
{Department of Mathematics \\ Temple University \\ \newline
1805 North Broad Street \\ Philadelphia, PA 19122}
\email{jacob.woods@temple.edu}

\date{August 11, 2023}


\begin{document}

\maketitle

\begin{abstract}
This paper develops low-dimensional moment methods for advective problems on networks of domains. The evolution of a density function is described by a linear advection-diffusion-reaction equation on each domain, combined via advective flux coupling across domains in the network graph. The PDEs' coefficients vary in time and across domains but they are fixed along each domain. As a result, the solution on each domain is frequently close to a Gaussian that moves, decays, and widens. For that reason, this work studies moment methods that track only three degrees of freedom per domain---in contrast to traditional PDE discretization methods that tend to require many more variables per domain. A simple ODE-based moment method is developed, as well as an asymptotic-preserving scheme. We apply the methodology to an application that models the life cycle of forest pests that undergo different life stages and developmental pathways. The model is calibrated for the spotted lanternfly, an invasive species present in the Eastern USA. We showcase that the moment method, despite its significant low-dimensionality, can successfully reproduce the prediction of the pest's establishment potential, compared to much higher-dimensional computational approaches.
\end{abstract}

\section{Introduction}
Numerous real-world processed can be described via network flows, consisting of multiple domains connected via flux coupling conditions, and the temporal evolution of density fields on the domains is described by a partial differential equation (PDE). Classical examples of such applications include traffic flow networks \cite{HoldenRisebro1995, HertyKlar2003, CocliteGaravelloPiccoli2005, GaravelloPiccoli2006, FarjounSeiboldMeshfree2013} or supply chain networks \cite{ArmbrusterMarthalerRinghofer2004, GoettlichHertyKlar2005, HertyKlarPiccoli2007}. Here, as a specific application, principled models for the life cycle of forest pests are considered. These pests develop through different life stages (the domains) and can traverse different developmental pathways (the coupling conditions across domains), see \S\ref{subsec:application_slf}. The temporal evolution\footnote{In computational mathematics, the term ``temporal evolution'' denotes how the solution of a PDE behaves in time, and this is how it is used herein. That is in contrast to biological literature, where ``evolution'' is commonly understood as the change in heritable characteristics over generations---an aspect not modeled here.} of the population density with respect to developmental age is described via an advection-diffusion-reaction equation on each domain, where the coefficients of the different terms vary across domains, and also vary with respect to time, but they are constant with respect to the developmental age. As described in \S\ref{subsec:mechanism}, this model structure renders the solution on each domain to frequently be close to a single Gaussian that widens and decays while it moves through the domain. While traditional PDE discretization methods such as finite volume schemes (cf.~\S\ref{subsec:other_methods}) can be applied to numerically approximate the model's governing equations, such approaches tend to require hundreds of degrees of freedom to accurately capture the solution's shape. Furthermore, for problems that exhibit strongly peaked age-distributions (i.e., very narrow Gaussians), the numerical resolution needed by traditional methodologies to achieve acceptable accuracy is particularly high.

This work considers an alternative methodology: moment methods that track only three degrees of freedom per domain, namely the lowest moments of the age-distribution on each domain. Mappings between these moments and the parameterization of a Gaussian are employed to approximate the coupled-PDE-model by a low-dimensional system that tracks only three degrees of freedom per domain. The basic application of this approach yields an ODE-based method in \S\ref{subsec:ode_scheme}. In addition, an asymptotic preserving method is developed in \S\ref{subsec:ap_scheme}, which can also handle strongly peaked solutions in a robust fashion.

The resulting moment method is then applied to a specific calibrated model for the life cycle of the spotted lanternly (SLF), introduced in \cite{lewkiewicz_temperature_2022}. The biologically relevant quantities of interest are used to assess the accuracy of the low-dimensional moment method (compared to a higher-dimensional reference method) in capturing the key dynamics across a wide range of model parameters. The SLF (\emph{lycorma delicatula}) is a species of planthoppers native to China. In 2014, it was introduced to Eastern Pennsylvania, and subsequently established and spread to several other northeastern US states. SLF can cause serious damage to plants and structures, severely impacting agriculture (grapes, apples, hops, etc.) and industry (hardwood, lumber, etc.). Due to these severe economic \cite{URBAN2019} and environmental \cite{KIM2011} concerns, modeling and prediction \cite{LADIN2023, JUNG2017}, as well as management \cite{LEE2019} and control \cite{PARK2009, CHOI2012} of SLF populations have become vital tasks.

This paper is organized as follows. In \S\ref{sec:problem_and_applications}, the structure of the governing equations and the ecology application are described. The equations that the chosen moments satisfy are then derived in \S\ref{sec:moments_properties_reconstruction}, and the reconstruction via Gaussians and the effective mapping between moments and reconstruction is described. These components then culminate in \S\ref{sec:moment_methods} into two versions of closed moment methods, including an asymptotic preserving version. This section also contrasts the moment methods against alternative numerical approaches. The performance of the moment methods is demonstrated in \S\ref{sec:results_test_problems}, via carefully designed test problems. Then, in \S\ref{sec:results_slf}, the moment method is applied to the calibrated SLF model, and maps of predictive establishment potential are produced. A closing discussion is given in \S\ref{sec:conclusions}, including an outlook on future research directions.

\vspace{1.5em}
\section{Mathematical Problem and Application}
\label{sec:problem_and_applications}
Here we describe the mathematical formulation of the considered problems, first in general terms in \S\ref{subsec:governing_equations}, and then for the specific ecology application in \S\ref{subsec:application_slf}.

\subsection{Governing equations}
\label{subsec:governing_equations}
We describe the network of coupled computational domains via the following mathematical framework. Let $\mathcal{G} = (\mV,\mE)$ be a finite directed graph with vertices $\mV$ and edges $\mE\subset \mV \times\mV$ such that if $(u,v)\in\mE$, then $(v,u)\not\in\mE$, i.e., between any pair of vertices, flow is allowed in one direction only. On each graph vertex $v \in \mV$, we consider a density function $\rho^v(a,t)$, defined on the domain $\Omega := (0,1)\times (0,\infty)$, that satisfies the governing linear advection-diffusion-reaction equation
\begin{equation}
\label{eqn:domain_PDE}
\pd{}{t}\rho^v + \pd{}{a}\!\prn{\nu^v(t)\rho^v - \xi^v(t)\pd{}{a}\rho^v} + \mu^v(t)\rho^v = 0\;,
\end{equation}
subject to the boundary conditions
\begin{equation}
\label{eqn:domain_boundary_conditions}
F^{v}_\text{in}(t) = f^{v}_\text{in}(t)\ \text{at}\ a=0\;,
\quad\text{and}\quad
\xi^v(t)\pd{}{a}\rho^v(1,t) = 0\;.
\end{equation}
Here $F_\text{in}(t) = \nu(t)\rho(0,t) - \xi(t)\pd{\rho}{a}(0,t)$ is the influx into the domain at $a = 0$ and it is assigned to equal $f_\text{in}(t)$, where $f_\text{in}$ is a given non-negative function of time. Note that in the diffusion-free case ($\xi^v = 0$), the boundary condition at $a = 1$ automatically becomes inactive.

For the first variable in $\rho^v(a,t)$ the letter $a$ is used, because in the SLF application (see \S\ref{subsec:application_slf}) it represents the developmental \emph{age}. Without loss of generality it is scaled to be $a\in [0,1]$ (in the application, $a$ represents the fraction of development achieved in a given life stage). The variable $t$ is time. In \eqref{eqn:domain_PDE}, the coefficients $\nu^v$, $\xi^v$, and $\mu^v$ are assumed non-negative at all times, and we further assume that on each domain, $\xi^v \le \mathcal{O}(\nu^v)$, i.e., as the advective speed $\nu$ approaches zero, the diffusion does so at the same rate (or faster).
In the SLF application, $\nu$ is the developmental speed, $\xi$ is a spreading of age distributions (see \S\ref{subsec:application_slf}), and $\mu$ is mortality. Suitable initial conditions $\rho^v(\cdot,0)$ are prescribed on all domains $v \in \mV$.

With the governing equations and boundary conditions defined on each domain, the domains/vertices are now coupled together (along the edges $\mE$) via fluxes. The outflux out of a domain (through $a = 1$) is
\begin{equation}
\label{eqn:domain_outflux}
F_\text{out}^v(t) = \nu^v(t)\rho^v(1,t) - \xi^v(t)\pd{}{a}\rho^v(1,t) = \nu^v(t)\rho^v(1,t)\;,
\end{equation}
where $(\xi\rho_a)$-term vanishes due the outflow boundary condition in \eqref{eqn:domain_boundary_conditions}. The coupling conditions are defined as follows. For a given vertex/domain $v \in \mV$, let
\begin{equation*}
I_v = \{u\in\mV \,:\, (u,v)\in\mE\}
\quad\text{and}\quad
O_v = \{w\in\mV \,:\, (v,w)\in\mE\}
\end{equation*}
denote the set of vertices that lead into $v$, respectively the set of vertices that flow out of $v$. If $I_v$ is empty, then $v$ is a source node, and an influx function $f_\text{in}$ is to be prescribed. For every edge $(v,w)\in\mE$, a ``split ratio'' or ``flux fraction'' $\alpha_{v,w}>0$ must be defined. If fluxes across domains exactly conserve mass, then one has for all $v \in \mV$ that $\sum_{w\in O_v}\alpha_{v,w} = 1$. However, it is also possible that not all the mass gets transferred across domains (e.g., mortality during life stage transitions), or that the mass multiplies (e.g., flux from egg-layers to eggs). Based on the split rations, the inflow into the domain $v$ is then given by
\begin{equation}
\label{eqn:coupling_fluxes}
f^v_\text{in}(t) = \sum_{u\in I_v}\alpha_{u,v} F_\text{out}^u(t)\;.
\end{equation}

A critical aspect is that the boundary/coupling conditions are formulated to preserve \emph{advective causality}, even if there exists some diffusive ``smearing'' in a domain. This is motivated by the application described in \S\ref{subsec:application_slf}: in the age development model, there is no backwards flow of information, that is: developmental age never decreases. The diffusive term models the empirical phenomenon that age distributions are observed to widen as development progresses, see \S\ref{subsec:application_slf}.

The above formulation avoids diffusive fluxes ($\xi \pd{}{a}\rho$) between domains. Moreover, because $\xi \le \mathcal{O}(\nu)$, there is no outflux if there is no aging. One should also note that the model, as first described in \cite{lewkiewicz_temperature_2022}, was formulated in terms of fractional steps: first, an advective sub-step is conducted, with coupled inflows and outflows, followed by a sub-step only pure diffusion and decay with zero flux at both boundaries. In the limit of the time step going to zero, the fractional steps and the model here are equivalent. Finally, the model in \cite{lewkiewicz_temperature_2022} also allows for an influx into a domain to result from the solution in another domain via a kernel, see below.

\subsection{Archetype mechanism of the network problem}
\label{subsec:mechanism}
A key building block for the network flows defined above is a single transition between two domains connected via a single edge $(v,w) \in \mE$, i.e., $O_v = \{w\}$ and $I_w = \{v\}$. We also assume $\alpha_{v,w} = 1$. For this scenario, let us call $v$ the \emph{sending domain} and $w$ the \emph{receiving domain}. Here the receiving influx equals the sending outflux, $F^w_\text{in} = F_\text{out}^v$. Furthermore, we consider the sending influx to be some prescribed function $f^v_\text{in}$, and the receiving outflux $F_\text{out}^w$ is irrelevant.

To understand the key mechanism caused by flux couplings at varying advection speeds, let us also consider the case of no diffusion, $\xi^v = 0 = \xi^w$, and no decay, $\mu^v = 0 = \mu^w$. Intuitively, this situation is like two conveyor belts transporting sand, with the first belt dropping sand onto the second, and sand grains tending to pile perfectly on top of each other. If the sending speed $\nu^v$ equals the receiving speed $\nu^w$, the density profile passes across the domain boundary unmodified. If the sending speed is slower than the receiving speed, i.e., $\nu^v < \nu^w$, the profile widens upon domain transition. Conversely, if the sending speed is larger, i.e., $\nu^v > \nu^w$, the profile sharpens, i.e., it compresses in the $a$-direction and increases in the $\rho$-direction.

In the limiting situation of zero receiving speed, i.e., $\nu^w = 0$ while $f^w_\text{in} > 0$, the influx boundary condition \eqref{eqn:domain_boundary_conditions}, which here would read as $F^w_\text{in} = \nu^w \rho^w(0,t) = f^w_\text{in}$, cannot be satisfied in the sense of functions. However, in many applications such a situation makes perfect practical sense, for instance: eggs keep getting produced, but no egg development occurs; or: incoming sand keeps piling up at the start of a stopped receiving conveyor belt. It is therefore reasonable to allow for solutions in the sense of distributions, i.e., allow for the creation of a Dirac delta distribution at $a = 0$. Clearly, the accurate representation of such solutions is numerically challenging; and one key advantage of moment methods is that they face no fundamental limitations in capturing this Dirac delta limit. In this work, we conduct a regularization by replacing Dirac peaks by narrow Gaussians (see \S\ref{subsec:moment_reconstruction}). However, we do allow these Gaussians to be extremely narrow. Hence, for applications in which there is a lack of widening diffusive effects, moment methods that are asymptotic preserving in this described regime are of critical importance, as described in \S\ref{subsec:ap_scheme}.

\subsection{Application: Life cycle of forest pests}
\label{subsec:application_slf}
For many species, in the modeling of its population dynamics it is advantageous to break the population into life-stages where the populations' dynamics (such as aging or death) are uniform within each life stage. For many insects such life stages come naturally due to their life cycle in instars, separated by molting events. Based on these life stages one can formulate a model for how a population may develop over time with a stage-structured system of age-structured PDE. This is a well-known methodology for modeling forest pest as well as other biological systems \cite{IANNELLI2017}. These systems often assume a form similar to equation \eqref{eqn:domain_PDE} but with varying flux coupling conditions between stages. Models of this kind have been calibrated for many forest pests \cite{GILIOLI2015, GILIOLI2017, PASQUALI2019}. In this work we focus on a specific model, calibrated for the spotted lanternfly \cite{lewkiewicz_temperature_2022}. In \S\ref{sec:results_slf}, this model is used as a test problem for the moment method developed in \S\ref{subsec:ode_scheme}, to assess the method's utility in simulating the life cycle of forest pests efficiently.

\begin{figure}
    \includegraphics[width=.90\linewidth]{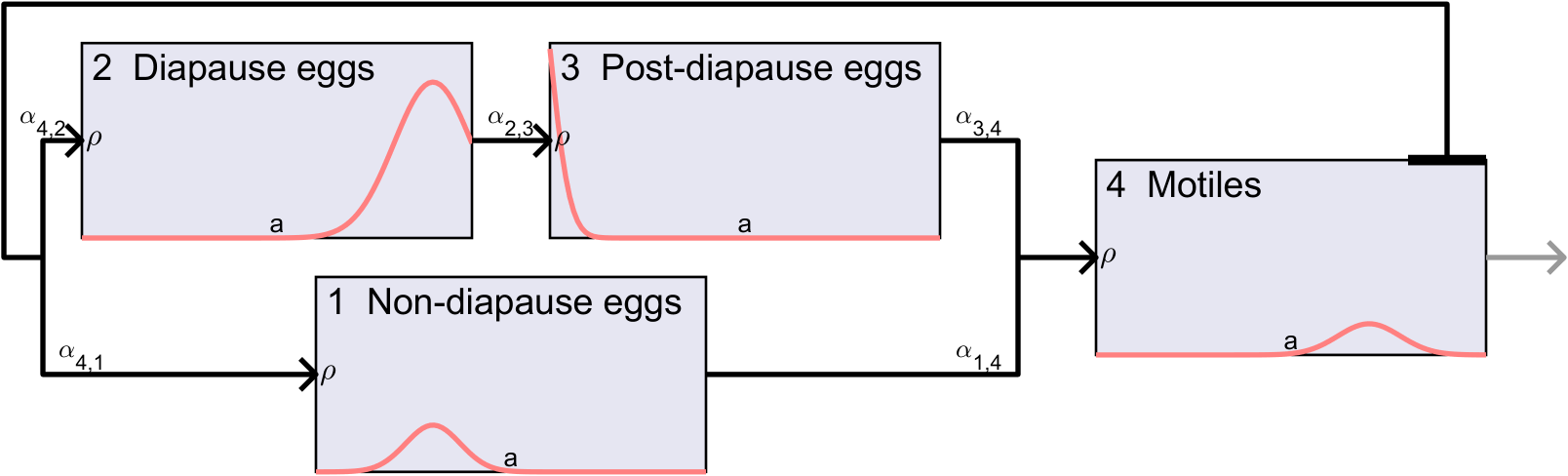}
    \caption{Example network of flux-coupled domains, here specific for the SLF life stage network used in \S\ref{subsec:application_slf}. Each of the four domains 1--4 has an age-density $\rho^v(a,t)$ described by a PDE of form \eqref{eqn:domain_PDE}, and out- and in-fluxes are distributed along arrows. For the SLF application, the flux from domain 4 to domains 1 and 2 is of a special structure, see equation~\eqref{eq:slf_model_egg_laying}. The gray arrow and the black arrow's base visualize this.}
    \label{fig:network_slf}
\end{figure}

In the model from \cite{lewkiewicz_temperature_2022}, the SLF life cycle is organized into four ($v \in \{1,2,3,4\}$) life-stage domains: three egg stages and one motile stage. The reason for this division is that in each domain, individuals share the same response to temperature which drives their development and mortality, i.e., the parameters $\nu$, $\xi$, and $\mu$ in \eqref{eqn:domain_PDE} are independent of the age variable $a$. However, the parameters do depend on time because they are functions of the ambient temperature $T$, i.e., for each domain there are functions $\nu^v(T)$, $\xi^v(T)$, and $\mu^v(T)$ (see \cite{lewkiewicz_temperature_2022}), and if a time-varying ambient temperature profile $T(t)$ is given, this effectively induces time-dependent coefficients in \eqref{eqn:domain_PDE}. The four stages, depicted in Fig.~\ref{fig:network_slf}, are:
\begin{enumerate}[\quad$\bullet$]
\item Domain 1: ``Non-diapause eggs''. Eggs develop normally, with the advection speed $\nu$ increasing with temperature (within certain limits \cite{lewkiewicz_temperature_2022}), and have normal resistance to extreme temperatures (expressed via the mortality function $\mu(T)$).
\item Domain 2: ``Diapause eggs''. Diapause is a state of slowed or paused development in exchange for improved cold resistance \cite{KOSTAL2006} in which warmth suppresses the advection speed $\nu$.
\item Domain 3: ``Post-diapause eggs''. Eggs' development is driven by temperature as in domain 1. However, an increased cold resistance from the diapause state is retained.
\item Domain 4: ``Motiles''. This combines all non-egg SLF: nymphs, instars, and adults, including egg-laying adults. At the motile stage, the history of the (non-)diapause pathway is lost, and the response to temperature is uniform.
\end{enumerate}
The transitions between life stages is modeled via the flux coupling conditions, depicted in in Fig.~\ref{fig:network_slf} via arrows. All depicted nonzero $\alpha_{v,w} = 1$, except for the arrows emanating from the motile motile domain: the outflow \eqref{eqn:domain_outflux} from domain 4 goes nowhere (death due to senescence); instead, eggs are produced via an egg laying kernel assigned to egg domains in a time-dependent fashion, see below.

On each domain $v \in \{1,2,3,4\}$, the function $\rho^v(a,t)$ represents the age density of females in that life stage. The stage-structured PDE assumes the form of \eqref{eqn:domain_PDE} with standard flux coupling boundary conditions for the fluxes out of domains 1, 2, and 3:
\begin{align*}
F_\text{in}^3 &= \nu^3(t)\rho^3(0,t) - \xi^3(t)\pd{}{a}\rho^3(0,t) = \nu^2(t)\rho^2(1,t) - \xi^2(t)\pd{}{a}\rho^2(1,t) = F_\text{out}^2\;, \\
F_\text{in}^4 &= \nu^4(t)\rho^4(0,t) - \xi^4(t)\pd{}{a}\rho^4(0,t)
= \nu^1(t)\rho^1(1,t) - \xi^1(t)\pd{}{a}\rho^1(1,t) \\
&\hspace{8em} + \nu^3(t)\rho^3(1,t) - \xi^3(t)\pd{}{a}\rho^3(1,t)
= F_\text{out}^1+F_\text{out}^3\;.
\end{align*}
The flux from domain 4 into domains 1 or 2 slightly deviates from the structure of \S\ref{subsec:governing_equations} as the production of eggs not be determined by the outflux of domain 4, but rather by an egg-laying kernel as follows.
\begin{equation}
\label{eq:slf_model_egg_laying}
\begin{split}
F_\text{eggs} &= \beta\nu^4(t)\int_{0}^{1} k(a) \rho^4(a,t) \ud{a}\;, \\
F_\text{in}^1 &= \nu^1(t)\rho^1(0,t) - \xi^1(t)\pd{}{a}\rho^1(0,t) = F_\text{eggs}I(t)\;, \\
F_\text{in}^2 &= \nu^2(t)\rho^2(0,t) - \xi^2(t)\pd{}{a}\rho^2(0,t) = F_\text{eggs}(1-I(t))\;.
\end{split}
\end{equation}
Here $I(t)$ is an indicator function that is $I(t) = 1$ over the duration of the year when eggs are laid into non-diapause and $I(t) = 0$ when eggs are laid into diapause. In \cite{lewkiewicz_temperature_2022}, $I(t)$ is modeled based on the photo-period of the day when the eggs was laid, with $I = 0$ after summer solstice but before winter solstice, and $I = 1$ otherwise. The kernel $k(a)$ describes the rate of egg production of a female at a given motile developmental age $a$. Its integral $\int_0^1 k(a) \ud{a}$ is the total number of eggs produced by a female if it lives until death by senescence. The constant $\beta\in [0,1]$ captures that there is some background mortality $1\!-\!\beta$ not related to temperature.

The assumption $\xi^v \le \mathcal{O}(\nu^v)$ introduced in \S\ref{subsec:governing_equations}, which justifies the ``causal'' flux couplings across domains, is satisfied in the SLF model as follows. The ratio $\xi^v/\nu^v$ may differ across domains, but it is always independent of temperature. This expresses the fact that age distributions of SLF observed in the field widen in a diffusive fashion, plausibly caused by the real-world fact that at different times, different individuals develop/age at different rates, even when exposed to the same circumstances. Here, it is important to emphasize that for cold-blooded species, developmental age is not identical to the time spent alive.

A key feature of the ecological model studied here is the presence of multiple pathways that individuals can take through the life stage network. The diapause pathway equips eggs with an increased resistance to cold, which can be crucial to enable enough eggs to survive the cold season in temperate climates. However, there is an additional ecological effect of diapause: it tends to synchronize the generational cycle with the annual temperature. Without diapause, eggs deposited in late summer will start developing and, given a warm enough fall, may hatch just went winter comes which will kill all hatched motiles. Conversely, diapause prevents the possibility of multiple generations per year, so in extremely warm climates, the non-diapause pathway may in turn be advantageous. Model simulations, as conducted in \S\ref{sec:results_slf}, are therefore critical towards understanding the life cycle mechanics and predicting the establishment potential of invasive species like SLF.

Ecological stage- and age-structured PDE models have been explored computationally prior, e.g., in \cite{HE2018, Pelovska2013RUNGEKUTTAMF, IANNELLI1997}, however, only via standard numerical methods (which track many numerical degrees of freedom per domain) and not exploring the significant model reduction offered by a well-designed moment method. Another important novelty is that the specific model from \cite{lewkiewicz_temperature_2022} provides pathways of co-existence of diapause and non-diapause that have not been widely explored numerically. The structure of the model, particularly the transitions between domains involving diapause, provides precisely the mechanism that tends to produce Gaussian density profiles, as described in \S\ref{subsec:mechanism}, at least in a wide range of relevant temperature profiles (see the results in \S\ref{sec:results_slf}).

\vspace{1.5em}
\section{Properties of Moments and Moments Reconstruction}
\label{sec:moments_properties_reconstruction}
This section establishes the building blocks used in \S\ref{sec:moment_methods} to formulate moment methods: the choice of moments and their governing equations are described in \S\ref{subsec:choice_of_moments}; and the needed transitions between moments and PDE solutions/reconstructions are established in \S\ref{subsec:moment_reconstruction}.

\subsection{Choice of moments and their governing equations}
\label{subsec:choice_of_moments}
Due to the structure of the governing equations described in \S\ref{subsec:governing_equations}, the considered models have a tendency to produce solutions that are close to Gaussians (that may decay and widen as they travel through the domains). Because a Gaussian is uniquely defined by three parameters, we choose moment methods that track exactly three moments on each domain in the network. Specifically, we track the lowest three monomial moments of the density function $\rho(a,t)$ on the domain $a\in (0,1)$, which are:
\begin{equation}
\label{eqn:monomial_moments}
m_0(t) = \int_{0}^{1} \rho(a,t) \ud{a}\;,\quad
m_1(t) = \int_{0}^{1} \rho(a,t)a \ud{a}\;,\quad\text{and}\quad
m_2(t) = \int_{0}^{1} \rho(a,t)a^2 \ud{a}\;.
\end{equation}
As the results for the SLF application in \S\ref{sec:results_slf} show, the assumption of the density on each domain being close to a Gaussian is indeed well satisfied in the current establishment region in the United States (see Fig.~\ref{fig:slf_results_Berks}), where the life cycle is dominated by diapause which synchronized species development and thus creates very peaked distribution functions $\rho$.

Closely associated with the monomial moments \eqref{eqn:monomial_moments} are the derived key quantities: (i) the mass $M$ of species in the domain; (ii) the mean (developmental age) $E$ in the domain; and (iii) the variance (of developmental age) $V$ in the domain, which are given in terms of the monomial moments as:
\begin{equation*}
M = m_0\;,\quad
E = \frac{m_1}{m_0}\;,\quad\text{and}\quad
V = \frac{m_2}{m_0}-E^2\;.
\end{equation*}
Computationally, we choose to track the monomial moments $(m_0,m_1,m_2)$, because their governing moment equations (see below) are easier. However, modulo the singularity at $M\to 0$, a description in terms of $(M,E,V)$ would be equivalently possible.

Numerical methodologies that represent the density function $\rho$ via pointwise values or cell averages (cf.~\S\ref{subsec:other_methods}) tend to require hundreds of degrees of freedom per domain. In comparison, tracking only three moments represents a significant reduction in dimensionality, which---if successful---can lead to substantial reductions in computational demand. In addition, there is a conceptual advantage of descriptions in terms of moments: the quantities total population $M$, mean developmental age $E$, and age variance $V$ are very intuitive and familiar concepts for ecologists and stakeholders working with these pests. While those moments could of course also be computed (post-hoc) for traditional methods like finite volumes, a reduced model and software that more directly track/update moments have the potential to be more accessible.

We now derive the equations that govern the rates of change of the moments \eqref{eqn:monomial_moments}. We differentiate \eqref{eqn:monomial_moments} with respect to $t$, invoke the governing PDE \eqref{eqn:domain_PDE}, employ integration-by-parts, and/or cancel/substitute quantities via the boundary conditions \eqref{eqn:domain_boundary_conditions}. For notational efficiency here we suppress domain indices and the time variable. In that abridged notation, \eqref{eqn:domain_boundary_conditions} reads as $(\nu\rho-\xi\rho_a)(0) = f_\text{in}$ and $\xi\rho_a(1) = 0$. With that, we obtain:
\begin{equation}
\label{eq:moments_evolution_equations}
\begin{split}
\frac{\ud{}}{\ud{t}}m_0 &= \int_0^1 \!\! -\nu\rho_a + \xi\rho_{aa} - \mu\rho \ud{a}
= -\nu[\rho]_0^1 + \xi[\rho_a]_0^1 - \mu m_0
= f_\text{in} - \nu\rho(1) - \mu m_0\;, \\
\frac{\ud{}}{\ud{t}}m_1 &= \int_0^1 \!\! -\nu a\rho_a + \xi a\rho_{aa} - \mu a\rho \ud{a}
= [-\nu a\rho + \xi a\rho_a - \xi\rho]_0^1 + \nu m_0 - \mu m_1 \\
&= -\nu \rho(1) - \xi (\rho(1) - \rho(0)) + \nu m_0 - \mu m_1\;, \\
\frac{\ud{}}{\ud{t}}m_2 &= \int_0^1 \!\! -\nu a^2\rho_a + \xi a^2\rho_{aa} \!-\! \mu a^2\rho \ud{a}
= [-\nu a^2\rho \!+\! \xi a^2\rho_a \!-\! 2\xi a\rho]_0^1 + 2\nu m_1 + 2\xi m_0 \!-\! \mu m_2\! \\
&= -\nu \rho(1) - 2\xi \rho(1) + 2\nu m_1 + 2\xi m_0 - \mu m_2\;.
\end{split}
\end{equation}
While some of the terms in the ODEs' right hand sides can be written in terms of the moments themselves, other terms require the knowledge of the actual solution $\rho$ of the PDE \eqref{eqn:domain_PDE}. As the goal is to formulate a numerical method that tracks (on each domain) solely the moments $(m_0,m_1,m_2)$, a reconstruction mapping is needed that assigns to a given vector of moments a reconstructed function. Such mappings are constructed in \S\ref{subsec:moment_reconstruction}. Then, in \S\ref{sec:moment_methods} the components are combined into two versions of closed moment methods that approximate the original network PDE problems.

\subsection{Moment reconstruction}
\label{subsec:moment_reconstruction}
Careful attention must be paid to the reconstruction of the Gaussian from the moments. 
The reconstruction mapping is necessary to evaluate the fluxes, and for many applications, the reconstructed population function is the quantity of interest. The construction herein is a special case of the Hausdorff moment problem \cite{mead_maximum_1984, talenti_recovering_1987}, with a particular attention to the fact that exactly 3 moments are considered. We wish to design methods to represent a Gaussian on each domain, parameterized by $C$, $a_0$, and $\sigma$:
\begin{equation}\label{eqn-gaussian-defn}
G(a;C,a_0,\sigma):=
       C\exp\prn{-\frac{(a-a_0)^2}{2\sigma^2}}.
\end{equation}
However, the moments $m_0$, $m_1$, and $m_2$ are being tracked. Hence, the mapping
\begin{equation*}
L : (C,a_0,\sigma) \mapsto (m_0,m_1,m_2)
\end{equation*}
and its inverse
\begin{equation*}
R : (m_0,m_1,m_2) \mapsto (C,a_0,\sigma)
\end{equation*}
must be understood and numerically approximated.
In particular, the methods constructed in \S\ref{sec:moment_methods} will need the reconstruction mapping
\begin{equation}
\label{eqn:Reconstruction}
\mR : (m_0,m_1,m_2) \mapsto G(\hspace{.5mm}\cdot\hspace{.5mm};C,a_0,\sigma)\;.
\end{equation}
The domain of this mapping (and others presented in this section) is a non-trivial issue, and will be discussed in \S\ref{subsubsec:reliazabilitydomain}. 

Due to the homogeneity of the mappings $L$ and $R$, it suffices to characterize the mappings between $(a_0,\sigma)$ and the statistical moments $(E,V)$, where $C$ may be chosen as a convenient non-negative constant:
\begin{equation*}
\LT : (a_0,\sigma) \mapsto (E,V)
\end{equation*}
and its inverse
\begin{equation*}
\RT : (E,V) \mapsto (a_0,\sigma)\;.
\end{equation*}
As neither the forward nor backwards map may be described with elementary functions, the mappings are implemented via a lookup table, moving the bulk of the computational complexity to a pre-computation that is conducted only once. Here we discuss some of the concerns that arise in the creation of the lookup table, and potential avenues to alleviate them.

\subsubsection{Realizability domain}
\label{subsubsec:reliazabilitydomain}
The first step is to characterize the \emph{realizability domain}, i.e., the set
\begin{equation*}
S = \{(m_0,m_1,m_2) : \exists\, C\ge 0, a_0\in\mathbb{R}, \sigma>0 \text{~s.t.~} (m_0,m_1,m_2) = L(C,a_0,\sigma)\}\;,
\end{equation*}
or, in the reduced variables,
\begin{equation*}
\tilde{S} = \{(E,V) : \exists\,  a_0\in\mathbb{R}, \sigma>0 \text{~s.t.~} (E,V) = \LT(a_0,\sigma)\}\;.
\end{equation*}
Below, we use $\mathbb{R}^+_0 := \mathbb{R}^+\cup \{0\}$. $S$ is defined as the image of $\mathbb{R}^+_0\times\mathbb{R}\times\mathbb{R}^+$ under $L$. Thus, the domain of $R$, the inverse mapping is $S$. Analogously, the domain of $\RT$ is $\ST$:
\begin{equation*}
L : \mathbb{R}^+_0\times\mathbb{R}\times\mathbb{R}^+ \rightarrow S
\quad\text{and}\quad
R : S \rightarrow \mathbb{R}^+_0\times\mathbb{R}\times\mathbb{R}^+ 
\end{equation*}
as well as
\begin{equation*}
\LT : \mathbb{R}\times\mathbb{R}^+ \rightarrow \ST
\quad\text{and}\quad
\RT : \ST \rightarrow \mathbb{R}\times\mathbb{R}^+ .
\end{equation*}

\begin{figure}
    \includegraphics[width=.70\linewidth]{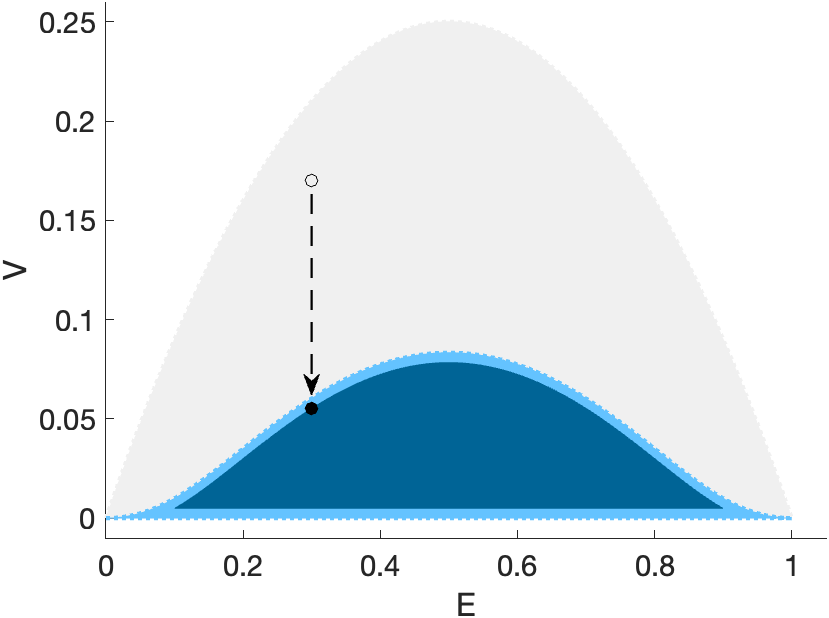}
    \caption{Realizability domains. The gray region represents all realizable statistical moments of positive functions on $(0,1)$. The light blue region is the realizability domain of Gaussians, $\ST$, and the dark blue region represents the regularized $\ST_D$. In this plot, the distance between $\ST$ and $\ST_D$ has been enlarged for clarity. A point outside $\ST_D$ is shown to be projected into the feasibility region $\ST$ to the boundary of $\ST_D$.}
    \label{fig:fullrealizability}
\end{figure}

We emphasize here that $S$ and $\ST$ refer to the realizability domain of the moments of single Gaussians on the unit interval. The realizability domain of the moments of arbitrary positive functions on the unit interval contains $S$ (see Fig.~\ref{fig:fullrealizability}). We recall that if $m_0> 0$,we can recover $E$ and $V$ from the moments. If $m_0 = 0$, the Gaussian is the constant zero function, and the values of $E$ and $V$ are not defined. It thus suffices to characterize the realizability domain of the statistical moments of the Gaussians. 

We then look to characterize the boundary of the set $\Tilde{S}$. As we consider positive functions restricted to $(0,1)$, for all $(E,V)\in\ST$, one has $E\in(0,1)$ and $V>0$. The curve that defines the top boundary, i.e., the value for which $V$ is maximized for a given $E$ value, cannot be represented as an elementary function $V(E)$. However, we can parameterize the curve, and then numerically determine the top boundary function, as follows.

\begin{figure}
    \includegraphics[width=.80\linewidth]{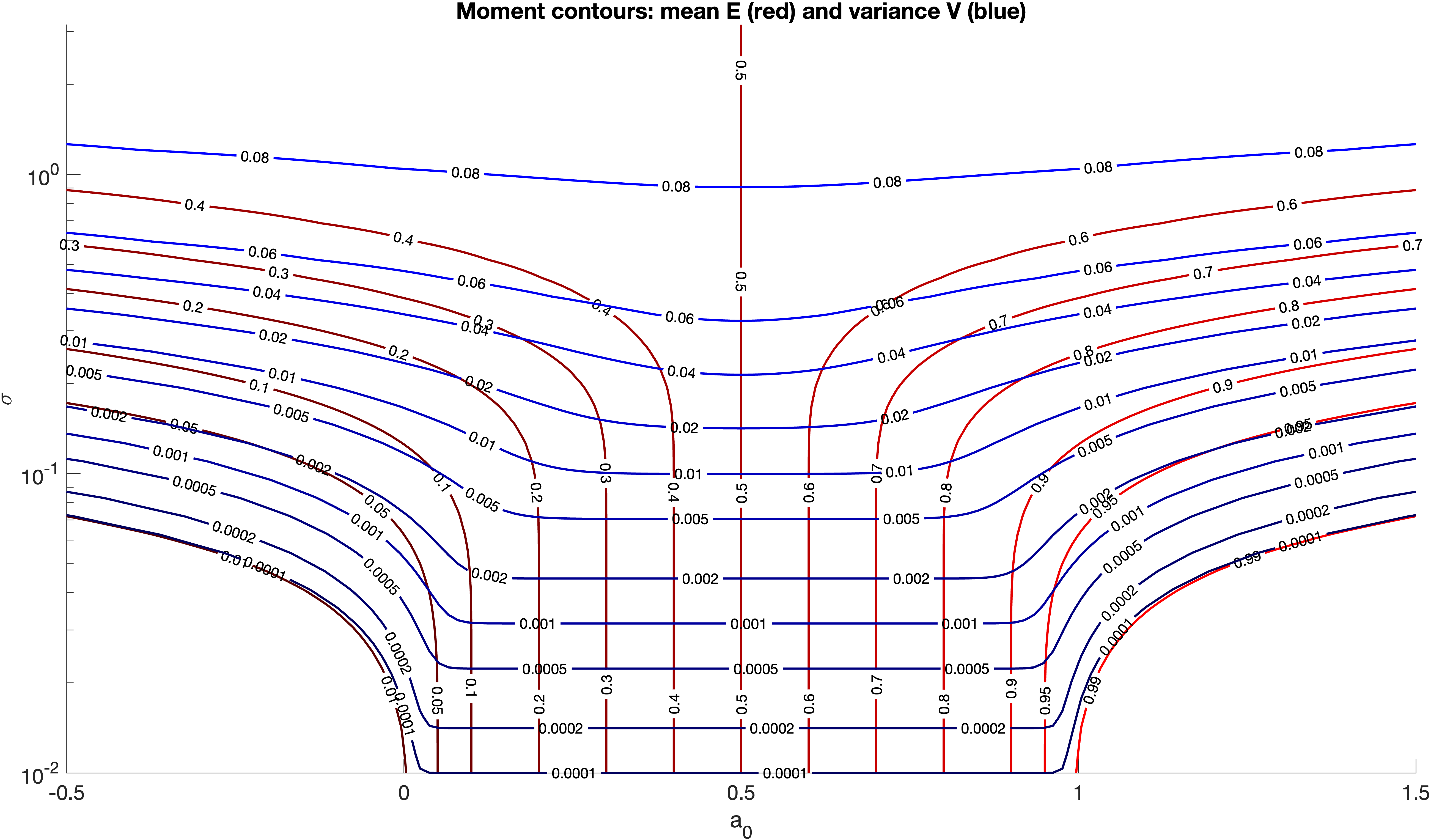}
    \caption{Contours of mean $E$ and variance $V$ of Gaussians on domain $[0,1]$, parameterized by the Gaussian peak location $a_0$ and width $\sigma$. Evidently, for a given mean $E$, the variance $V$ is maximized as $|a_0|\to\infty$.}
    \label{fig:contours_E_V}
\end{figure}

Figure~\ref{fig:contours_E_V} displays the $E$- and $V$-contours in the $(a_0,\sigma)$-domain. The plot illustrates many important structural properties. First, intuitively, for any fixed $a_0$, in the limit $\sigma\to\infty$, the Gaussian approaches a constant function and thus $(E,V)\to (\frac{1}{2},\frac{1}{12}$). Next, the plot illustrates the fact that, for a given mean $E$, the variance $V$ is maximized as $|a_0|\to\infty$. Consider WLOG the case $a_0>\frac{1}{2}$. Then, $E>\frac{1}{2}$, and the variance-maximizing Gaussians are obtained in the limit $a_0\to \infty$. In that asymptotic limit, we have on the domain $[0,1]$ that
\begin{equation*}
\exp\prn{-\frac{(a-a_0)^2}{2\sigma^2}} \cong \exp\prn{-\frac{(a_0)^2}{2\sigma^2}} \exp\prn{-\beta a}\;,
\end{equation*}
where $\beta = \frac{-a_0}{\sigma^2}$. The constant $\exp (-(a_0)^2/2\sigma^2)$ is irrelevant for the quantities $E$ and $V$. Due to symmetry, the case $a_0<\frac{1}{2}$ is treated analogously, hence the variance-maximizing distributions are simple exponential functions $\exp(-\beta a)$, parameterized by $\beta\in\mathbb{R}$.
Straightforward calculations yields that for $\beta\neq 0$,
\begin{equation*}
m_0(\beta) = \frac{1-e^{-\beta}}{\beta}
\;,\quad
m_1(\beta) = \frac{1-(1+\beta)e^{-\beta}}{\beta^2}
\;,\quad
m_2(\beta) = 2\frac{1-(1+\beta+\frac{1}{2}\beta^2)e^{-\beta}}{\beta^3}\;,
\end{equation*}
and $m_0(0) = 1$, $m_1(0) = \frac{1}{2}$, and $m_2(0) = \frac{1}{3}$. This yields the upper boundary of the realizability domain in the $(E,V)$-plane as the curve
\begin{equation*}
(E(\beta),V(\beta)) = \prn{\frac{m_1(\beta)}{m_0(\beta)},\frac{m_2(\beta)}{m_0(\beta)}-E(\beta)^2}\;.
\end{equation*}
In turn, the bottom boundary $(E,V) = (E,0)$ is formed by the distributions that for each $E\in [0,1]$ have minimal variance. Those would be Dirac deltas $\delta(a-a_0)$, which are the limits of Gaussians as $\sigma\to 0$. 

When computing the forward map for the lookup table, we select a finite number of points from the set of all possible $(a_0,\sigma)$ points, i.e., the open upper half plane. We seek to truncate the allowable $(a_0,\sigma)$ pairs without significantly reducing the size of the $\ST$. Let $D\subset\mathbb{R}\times\mathbb{R^+}$ be closed and bounded and 
\begin{equation}
\tilde{S}_D := \{(E,V) : \exists\,  (a_0,\sigma)\in D\text{~s.t.~} (E,V) = \LT(a_0,\sigma)\}.
\end{equation}
We seek $D$ such that $\ST_D$ is close to $S$.




\subsubsection{The forward mapping} \label{subsubsec:forwardmap}
Computing the forward mapping requires truncating allowable $(a_0,\sigma)$ pairs. However, our choice of pairs must be sufficiently rich to yield sufficient accuracy in the backwards mapping. Additionally, the mapping is prone to numerical errors if implemented na\"{i}vely. We will discuss possible approaches that will help to mitigate computational expenses and numerical errors without sacrificing the accuracy of the inverse lookup table. 

First, we leverage symmetry of the Gaussians and the domain about the point $(\frac{1}{2},0)$ and perform calculations only for Gaussians centered above this point, $a_0\geq\frac{1}{2}$. Symmetry gives that if $a_0'=1-a_0,$ and $\sigma'=\sigma$, the corresponding moments become $E'=1-E$ and $V'=V$. 

Numerically calculating the integral of a Gaussian over a finite interval accurately requires caution. The simplest approach is to use explicit formulas in terms of the Gauss error function: 
\begin{align}
\label{eqn:Gauss_m0}
m_0(C,a_0,\sigma) &= C\sigma\sqrt{\frac{\pi}{2}}\prn{\erf{\frac{1-a_0}{\sqrt{2}\sigma}}-\erf{\frac{-a_0}{\sqrt{2}\sigma}} }, \\
\label{eqn:Gauss_m1}
m_1(C,a_0,\sigma) &= a_0m_0-\sigma^2\prn{G(1;C,a_0,\sigma)-G(0;C,a_0,\sigma)}\;, \\
\label{eqn:Gauss_m2}
m_2(C,a_0,\sigma) &= a_0m_1-\sigma^2m_0-\sigma^2G(1;C,a_0,\sigma)\;.
\end{align}
This approach, though simple, may lead to significant numerical roundoff errors: when $\abs{x} > 6$, then $|\erf{x}| = 1$ within machine precision. Hence, if $a_0 \ll 0$ or $a_0 \gg 1$ while $\sigma \ll 1$, directly evaluating \eqref{eqn:Gauss_m0}--\eqref{eqn:Gauss_m2} will lead to unacceptable errors.

We may circumvent using the error function in favor of Simpson's rule for numerical integration. However, computing the integral of a Gaussian centered in the unit interval with $\sigma \ll 1$ requires a huge number of quadrature points to resolve. In this case, efficiency can be retained either by restricting the integration domain to where the Gaussian is greater than a desired tolerance (e.g., machine precision), or by replacing the integral over the finite domain to an integral of the Gaussian over all real numbers.

\begin{remark}
Computation of the forward map during the creation of the lookup table is considered a pre-computation, thus efficiency is not a concern. If the recovery of the mass is required as a runtime calculation, we may wish to implement a hybrid $\erf$-quadrature rule method to ensure efficiency, taking the error function when the arguments are sufficiently small, and computing the interval with quadrature when not. 
\end{remark}

\begin{figure}
    \centering
    \begin{subfigure}
        \centering
        \includegraphics[height=2.1in]{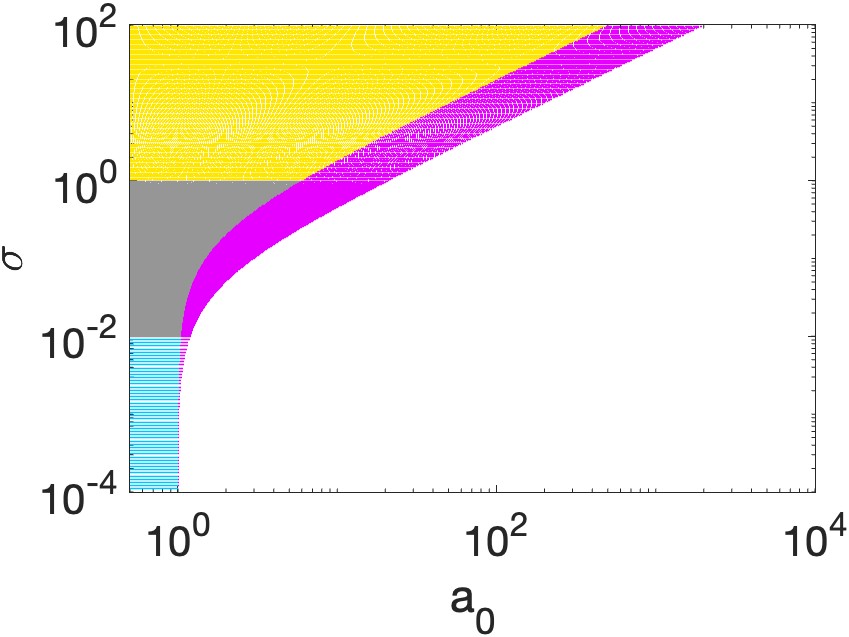}
    \end{subfigure}%
    \begin{subfigure}
        \centering
        \includegraphics[height=2in]{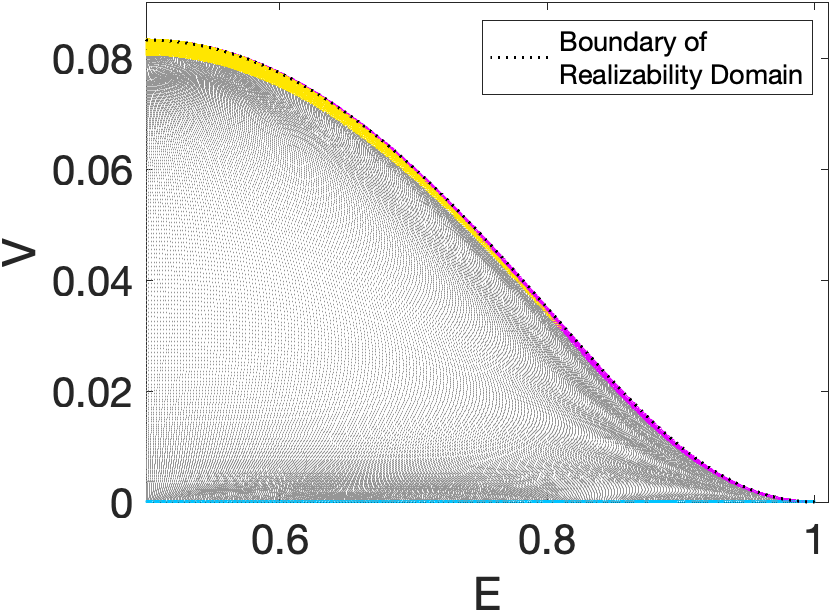}
    \end{subfigure}
    \caption{Points $(a_0,\sigma)\in D$ and the corresponding points in $(E,V)$. Cyan points represent functions close to Dirac delta distributions. Magenta points represent Gaussians centered far from $(0,1)$ with respect to their width (large values of $\beta$). Yellow points represent Gaussians with large width comparative to $a_0$ (small values of $\beta$). The gray points which make up the majority of the feasibility region in $(E,V)$ represent Gaussians of moderate width not far from the interval $(0,1)$.}\label{fig:a0sigma}
\end{figure}

Finally, we use unevenly spaced $(a_0,\sigma)$ points to more efficiently capture the realizability domain (see Fig.~\ref{fig:a0sigma}). We require that $\sigma$ is bounded away from $0$. In implementation, we vary $\sigma$ between $10^{-4}$ and $10^2$, with a higher concentration of values between $10^{-2}$ and $1$. For any given $\sigma$, we choose $a_0$ between $\frac{1}{2}$ and $1+c\sigma$. This uneven grid spacing avoids the unnecessary calculation of points when $a_0$ is very far from $[0,1]$ while $\sigma \ll 1$, and but also allows for $(a_0,\sigma)$ pairs with $a_0$ far from $[0,1]$ while $\sigma \gg 1$. We choose $c=20$ to balance computational cost and capturing the realizability domain well.

\subsubsection{The backwards map}
Mapping from the moments to the restriction of a Gaussian to the domain $(0,1)$ is ill-conditioned near the boundary of the realizability domain. In particular, though any $(E,V)$ pair corresponds uniquely to a Gaussian reconstruction, we may find $(E,V)$ and $(E',V')$ that are arbitrarily close, while the corresponding Gaussian parameters $(a_0,\sigma) = \RT(E,V)$ and $(a_0',\sigma') = \RT(E',V')$ are arbitrarily far from each other. In light of this, accurate reconstruction may seem hopeless.
However, though $a_0$ and $\sigma$ are extremely sensitive to perturbations in $(E,V)$, we do not directly use these values in our moment method, but rather the reconstructed Gaussian. The pointwise values of the reconstructed Gaussians are in general not as sensitive to small perturbations as $a_0$ and $\sigma$. 

A simple example is to consider points near the bottom boundary: for $\sigma$ small, $G(a;C,a_0,\sigma)$ approximates a Dirac delta on $(0,1)$, i.e., it is near zero a small distance from $a_0$. Similarly, we can take $G(a;C,a_0,c\sigma)$, with $c$ large (but so $c\sigma\ll 1$.) Then, outside the small neighborhood around $a_0$, there is near perfect agreement of the pointwise values of the functions. 

\begin{remark}
Though the pointwise values of reconstructed Gaussians are not as sensitive to perturbations as the values of $a_0$ and $\sigma,$ care must be taken when using reconstructed values. If $a_0$ is close to 1, small perturbations in $\sigma$ will have a large effect on the evaluation of the reconstructed Gaussian at $a=1$, which may effect the accuracy of the moment method \eqref{eq:moment_system_ode}.
\end{remark}

\begin{figure}
    \centering
    \begin{subfigure}
        \centering
        \includegraphics[height=2in]{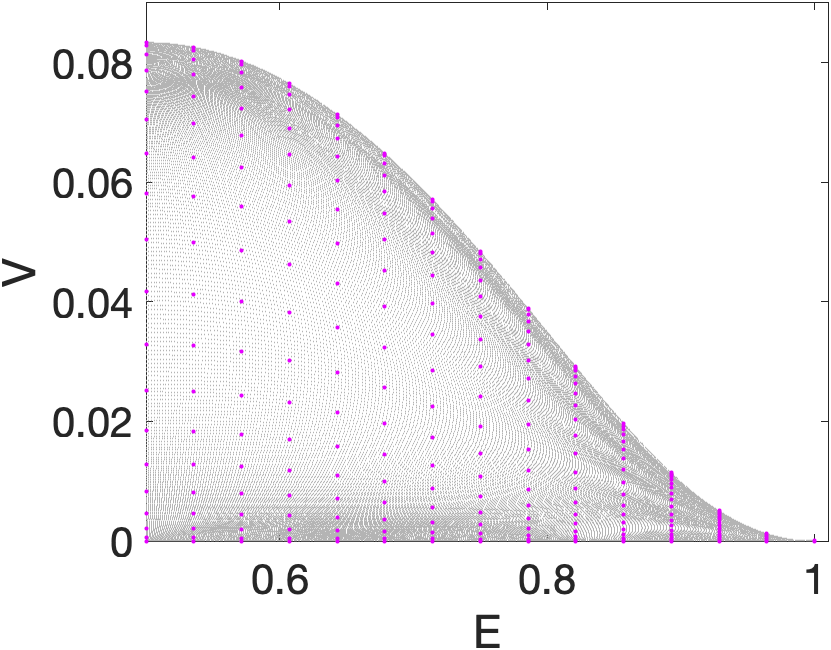}
    \end{subfigure}%
    \caption{Points selected in the $(E,V)$ are seen to be concentrated near the boundary. The displayed points are a selection of the total points, thinned for clarity.}\label{fig:backwardsmapEV}
\end{figure}

For the points stored in the lookup table, a regular grid in $(E,V)$ would be inefficient, as we need high resolution along the boundary of the realizability domain. Instead, we create a fine, uniform grid in $E$, and for each grid point $e_i$, we create an irregular grid in $V$ from $v=0$. to $v=\max_{(e_i,v^*) \in \Tilde{S}}v^*$ with points being more concentrated near the top and bottom boundary (see Fig.~\ref{fig:backwardsmapEV}). 

After selecting points, we create the lookup table by linear interpolation of the points generated by the forward map. The full realizability domain of $(E,V)$ cannot be captured by a finite number of points. Thus, if a point lies outside of the range of our forward mapping, we implement nearest neighbor interpolation. (This is equivalent to the projection described in \ref{subsubsec:regularization}.) We set the constant
\begin{equation*}
    C^*:=\int_0^1G(a;1,a_0,\sigma)\ud{a}\;,
\end{equation*}
and calculate it directly from the recovered $(a_0,\sigma)$ values. We note that it is also possible to calculate these constants with the initial forward map and to recover $C^*$ via interpolation. However, then interpolation errors might lead to a $(a_0,\sigma,C^*)$ triple such that $C^*\neq\int_0^1G(a;a_0,\sigma)\ud{a},$ which in turn may cause errors particularly in the asymptotic preserving method (see \S\ref{subsec:ap_scheme}).

\begin{remark}\label{remark:recoveryscheme}
A simple recovery scheme using the lookup table is as follows. Given $(M,E,V)$, we recover the corresponding Gaussian $(a_0,\sigma,C)$. 
\begin{enumerate}
    \item Find $e_i$, the nearest neighbor to $E$ (one-dimensional nearest neighbor interpolation).
    \item Find $v_{i,j}$, the nearest neighbor to $V$ along the $i^{th}$ row (one dimensional nearest neighbor interpolation).
    \item Recover $(a_0)_{i,j},\sigma_{i,j},C^*_{i,j}$ from the lookup table.
    \item Set $G(a;a_0,\sigma,C)=\frac{m_0}{C_{i,j}^*}\exp\prn{-\frac{(a-(a_0)_{i,j})^2}{2\sigma^2_{i,j}}}$. 
\end{enumerate}
This approach ensures that $m_0 = \int_0^1 G(a;a_0,\sigma,C)\ud{a}$ (critical for mass conservation), while being less costly than more complicated approaches like interpolation on the uneven grid. It balances computational cost and accuracy, and allows for sizable lookup tables with significant refinement near the boundary. 
\end{remark}

\subsubsection{Regularization and extension}
\label{subsubsec:regularization}
A naturally question is: how to address a given set of moments that are not in the feasible domain? This may occur due to numerical errors (e.g., from a time-stepping scheme used), or for systematic reasons, e.g., the true solution may not be Gaussian and its moments may not be representable by a Gaussian. One simple scenario for the latter case is a time-dependent advection speed $\nu(t)$ in a domain that receives a constant influx. If $\nu$ starts slowly, then speeds up, and then slows down again before the mass has left the domain, the solution will naturally form a bimodal distribution whose variance could easily larger than what is realizable via Gaussians.

We recall the definition of $\ST_D$ and further define the set of feasible moments of Gaussians determined by $(a_0,\sigma) \in D$ as
\begin{equation}
S_D = \{(m_0,m_1,m_2) : \exists\, (a_0,\sigma)\in D,\ C\in \mathbb{R}^+\text{~s.t.~} (m_0,m_1,m_2)= L(C,a_0,\sigma)\}\;.
\end{equation}

Note that if $(E,V)\in\ST_D$ and $m_0\geq 0$, then $\prn{m_0,m_0 E,m_0 (V+E^2)}\in S_D$. thus this space is essentially equivalent to $\ST_D$. Boundary values of $\ST_D$ are defined as $E_\text{min} = \min_{E\in \ST_D}E$ and $E_\text{max} = \max_{E\in \ST_D}E$. The minimal and maximal $V$ values are dependent on $E$, thus we define $V_{\text{min},E}=\min_{(E,V)\in\ST_D}V$ and $V_{\text{max},E}=\max_{(E,V)\in\ST_D}V$. With this, we define the projection $\Tilde{P}:\mathbb{R}\times\mathbb{R}\rightarrow\ST$ as
\begin{equation*}
\Tilde{P}(E,V) =
\begin{cases}
    (E,V) & (E,V)\in\ST_D \\
    (E_\text{min},V_{\text{min},E_\text{min}}) & E<E_\text{min} \\
    (E_\text{max},V_{\text{max},E_\text{max}})\quad & E>E_\text{max} \\
    (E,V_{\text{min},E}) & V<V_\text{min} \\
    (E,V_{\text{max},E}) & V>V_\text{max}\;.
\end{cases}
\end{equation*}
We note that $E_\text{min}$ and $E_\text{max}$ are close to 0 and 1 respectively, and that if $\rho(a)\geq 0$ on $(0,1)$, $E$ should not leave the interval $(0,1)$ (numerical approximation errors notwithstanding). However, the $\ST$ is contained in the set of feasible statistical moments for all positive functions, thus $V$ may achieve values outside the feasibility region without any numerical errors occurring. In this case, we project down towards the boundary of $\ST_D$, leading to a function with maximal possible variance (we recall $a_0$ and $\sigma$ both approach infinity at the boundary of $\ST$).
By projecting negative and near zero values for $V$ up to $\ST_D$, Dirac delta distributions are regularized and replaced with narrow Gaussians.   
This projection allows the method to be defined even when the moments are not in the feasibility region $S$. 

Let $(m_0,m_1,m_2)\in \mathbb{R}^3$. First, if $m_0\leq0,$ we are not in the feasibility region, and we set $(m_0,m_1,m_2) = (0,0,0)$, corresponding to the zero function (even though $E$ and $V$ are not defined in this case). Next, if $m_0>0$, we calculate $(E_m,V_m)=\prn{\frac{m_1}{m_0},\frac{m_2}{m_0}-\prn{\frac{m_1}{m_0}}^2}$, project onto the feasible region $\ST_D,$ and from there recover the moments
\begin{equation}
\label{eq:projected moments}
(m_0,m_1,m_2) =
\begin{cases}
    m_0(1,\PT E_m,\PT V_m+(\PT E_m)^2)\quad & m_0>0 \\
    (0,0,0) & m_0\leq 0\;.
\end{cases}
\end{equation}
This allows us to define a backwards mapping $\bar{R}:=RP$ from any triple, realizable or not, to the Gaussian parameters,
\begin{equation*}
\bar{R} : \mathbb{R}^3 \rightarrow \mathbb{R}^+_0 \times \mathbb{R}\times\mathbb{R}^+\;,
\end{equation*}
with
\begin{equation}\label{eqn:Reconstruction_regularized}
\bar{\mathcal{R}}(m_0,m_1,m_2)(a) = \frac{m_0}{C^*}\exp{\prn{-\frac{(a-a_0)^2}{2\sigma^2}}}\;.
\end{equation}
For ease of notification, in future sections we shall suppress the bar notation, and refer to $\mathcal{R}(m_0,m_1,m_2)(a)$ as the Gaussian recovered through the projection and lookup table evaluated at $a,$ and $R(m_0,m_1,m_2)$ as the parameters of the Gaussian.

\vspace{1.5em}
\section{Moment Methods}
\label{sec:moment_methods}
Based on the building blocks established in the prior sections, we now formulate two versions of closed moment methods. We construct an ODE-based moment method in \S\ref{subsec:ode_scheme} that is conceptually simple and works well if the reconstructed functions are not overly peaked. Then, in \S\ref{subsec:ap_scheme} we present a more specialized scheme that is asymptotic preserving in that it works in a robust fashion with large time step sizes even for very peaked solutions. The role of these moment methods in comparison with other numerical methods is discussed in \S\ref{subsec:other_methods}.

\subsection{ODE-based moment method}
\label{subsec:ode_scheme}
We now combine the governing equations for the moments \eqref{eq:moments_evolution_equations} with the (regularized) reconstruction mapping \eqref{eqn:Reconstruction} to obtain the system of ODEs, in which for each domain $v\in \mV$ one has
\begin{equation}
\label{eq:moment_system_ode}
\begin{split}
\tfrac{\ud{}}{\ud{t}}m_0^v &= - \nu^v \rho^v(1) - \mu^v m_0^v + f_\text{in}^v\;, \\
\tfrac{\ud{}}{\ud{t}}m_1^v &= -(\nu^v+\xi^v) \rho^v(1) + \xi^v \rho^v(0) + \nu^v m_0^v - \mu^v m_1^v\;, \\
\tfrac{\ud{}}{\ud{t}}m_2^v &= -(\nu^v+2\xi^v) \rho^v(1) + 2\nu^v m_1^v + 2\xi^v m_0^v - \mu^v m_2^v\;.
\end{split}
\end{equation}
where we use the short notation for the reconstruction $\rho^v(a) = \mathcal{R}(m_0^v,m_1^v,m_2^v)(a)$. Moreover, via \eqref{eqn:domain_outflux} the outflux from domain $v$ is $F_\text{out}^v = \nu^v \rho^v(1)$, again using the reconstruction. The outfluxes from all domains, via the flux coupling conditions \eqref{eqn:coupling_fluxes}, determine the influxes $f_\text{in}^w$ for all domains $w\in \mV$.

This yields a $3n$-dimensional ODE system, where $n$ is the number of domains. Within each domain, the three moments are coupled via the reconstruction mapping $\mathcal{R}$, and the equations across domains couple via \eqref{eqn:coupling_fluxes} from outfluxes to influxes. An important fact is that the ODE system \eqref{eq:moment_system_ode} is nonlinear, even though the underlying network PDE model is linear (see the discussion in \S\ref{subsec:other_methods}). Due to the regularized reconstruction \eqref{eqn:Reconstruction_regularized}, the right hand side of \eqref{eq:moment_system_ode} is always well-defined. However, it is not a-priori guaranteed that the solution of \eqref{eq:moment_system_ode} remain close to the true moments of the solution of the PDE model. The numerical results in \S\ref{sec:results_test_problems} and \S\ref{sec:results_slf} demonstrate the success and possible modes of failure of the moment method.

It is also possible that a given initial condition for the PDE \eqref{eqn:domain_PDE} is not Gaussian, or even not realizable. Hence, we project the moments of any initial condition via \eqref{eq:projected moments} before starting the computation.

A key structural advantage of the moment method approach presented herein is that it generates a generic ODE system which in principle can be advanced via any time stepping method of choice. In particular, existing methodologies and toolboxes for high-order-in-time, adaptive time step choice, dense output, or even implicit time stepping \cite{WannerHairer1991} could easily be employed. Moreover, the structure of \eqref{eq:moment_system_ode} can even be amenable to positivity-preserving ODE solvers \cite{BurchardDeleersnijderMeister2003}. In turn, if no such specialized ODE solver is used, time-stepping errors could produce non-realizable moments. Of particular concern is the possibility that $m_0$ could become negative in a domain. In Runge-Kutta schemes, the constraint $m_0\ge 0$ can be ensured via an ad-hoc flux limiting applied in each stage: any outflux \eqref{eqn:domain_outflux} is capped to $F_\mathrm{out}^v = \min(\nu^v(t)\rho^v(1,t),m_0^v(t)/\Delta t)$, where $\Delta t$ is the time step. This ensures that the mass that leaves a domain never exceeds the mass that is inside the domain.

In all numerical tests conducted in \S\ref{sec:results_test_problems} and \S\ref{sec:results_slf}, the use of a simple RK4 with uniform time steps and the above flux limiting generated positive solutions; moreover, for reasonably small time steps, it even retained realizability, even though that property is not guaranteed. That is in contrast to automated methods like Matlab's \texttt{ode45.m} without flux limiting. For sharply peaked Gaussians, negative $m_0$-values frequently emerged; and even when that did not happen, the adaptive time stepper sometimes terminated without success because the target tolerance could not be met. For the same peaked Gaussian problems, RK4 with flux limiting always generated reasonable solutions, albeit not always at a desirable accuracy, see the results in Fig.~\ref{fig:APvsMOL}.



With only three degrees of freedom per domain, the dimensionality of the ODE system \eqref{eq:moment_system_ode} is likely as low as a reasonable reduced model of the full network PDE model can be, and the accuracy of this description is enabled by the problem-specific mechanism described in \S\ref{subsec:mechanism}. That being said, low-dimensionality does not automatically translate into low computational cost. One right hand side evaluation of \eqref{eq:moment_system_ode} requires (at least) one reconstruction per domain, whose fast evaluation relies on an efficient lookup table, created in \S\ref{subsec:moment_reconstruction}.

\subsection{Asymptotic preserving scheme for advection of narrow Gaussians}
\label{subsec:ap_scheme}
As described in \S\ref{subsec:mechanism}, the network advection problems studied herein possess a natural mechanism to create strongly peaked distributions. These pose computational challenges for the ODE-based moment method from \S\ref{subsec:ode_scheme}. Here we develop a numerical methodology that remedies some of these challenges. Since the problem is fully rooted in the advection part of \eqref{eqn:domain_PDE}, in this section we restrict the formulation to \eqref{eqn:domain_PDE} without diffusion or decay, i.e., $\xi$ and $\mu$ are zero. 

Consider the case of a narrow Gaussian of width $\sigma \ll 1$ traveling in a domain of speed $\nu = \mathcal{O}(1)$ to the domain's outflow boundary. This induces the characteristic time scale $\tau = \sigma/\nu$, on which solution values change at fixed $a$-locations. In order to accurately resolve this time scale, non-specialized time-stepping methods for the ODE-based moment methods tend to require time step sizes $\Delta t \leq \tau$. Hence, when $\tau$ is very small, as is the case when $\sigma \ll 1$ and $\nu = \mathcal{O}(1)$, traditional time-stepping methods become computationally expensive. A desirable numerical scheme tracks very narrow Gaussians without incurring the above restriction. Such numerical schemes that also work well when $\Delta t \gg \tau$, but sufficiently small for any accuracy requirements, are denoted \emph{asymptotic preserving} \cite{jin1999AP}.
The asymptotic preserving property is critical particularly in the case when flow of mass in a network of domains where the advection speed $\nu^v(t)$ in each domain $v$ can vary, creating Dirac delta distributions that must be accurately tracked. The AP scheme here tracks the mass that is ``pushed'' from one domain to another during one time step, effectively shifting the temporal rates of change (as in the the ODE-based method in \S\ref{subsec:ode_scheme}) to rates of change in the $a$-variable that can be resolved via direct formulas or simple quadratures.

To derive the method, let $\rho^v(a,t)$ be the density function in domain $v\in \mV$ (considering at least two domains) for $a\in [0,1]$, time $t$, and traveling with speed $\nu^v(t)$. Given corresponding moments $(\mzero^v(t),\mone^v(t),\mtwo^v(t))$ at time $t$, as defined in \eqref{eqn:monomial_moments}, we wish to determine the value of the moments at time $t+\Delta t$, accounting for the mass flow across domains. A key quantity is the length scale over which information is advected in domain $v$ from $t$ to $t+\Delta t$, which is:
\begin{equation*}
   \dela^v = \int_{t}^{t + \Delta t} \nu^v(t) \ud{t}.
\end{equation*}
For general unknown velocity functions $\nu^v(t)$, the integral needs to be approximated numerically. A simple way is to use $\dela^v \approx \Delta t\, \nu^v(t + \Delta t/2)$. However, more accurate quadratures are also possible. The mass outflux from domain $v$ is given by $\int_{1-\Delta a^v}^{1} \rho^v(a,t) \ \ud{a}$ where \begin{equation*}
\rho^v(a,t) = \mathcal{R}(\mzero^v(t),\mone^v(t), \mtwo^v(t)(a)\;.
\end{equation*}
For computational accuracy and efficiency, we equivalently write the mass outflux as
\begin{equation*}
\Mzero^v = \int_0^{\dela^v} \rho^v(1-\dela^v + a,t) \ud{a}\;.
\end{equation*}
Using the same shift and change of variables we additionally define the following convenient formulations for the first and second moment respectively:
\begin{equation*}
\Mone^v = \int_0^{\dela^v} a \rho^v(1-\dela^v + a,t) \ud{a}\;,
\quad
\Mtwo^v = \int_0^{\dela^v} (a)^2 \rho^v(1-\dela^v + a,t) \ud{a}\;.
\end{equation*}
The implementation of these integrals is crucial to the working of the numerical method. Note that this is the same problem addressed in \S\ref{subsubsec:forwardmap} on computing the forward map which can be done via analytical Gaussian error function formulas or numerical quadrature. The computation of the integrals is pseudo-exact, incurring numerical errors only when numerical quadrature is used. 

We now describe the changes due to the flow of mass between a sending domain $v$ and receiving domain $w$ that are part of an edge $(v,w) \in \mE$. To obtain the correct behavior in domain $w$, when domain $v$ and domain $w$ have different time-dependent speeds, we define the ratio $\displaystyle \gamma^{v \to w} = \frac{\nu^w(t + \Delta t/2)}{\nu^v(t + \Delta t/2)}$ (which approximates $\displaystyle \frac{\nu^w}{\nu^v}$ on the interval $[t,t+\Delta t]$ with $\mathcal{O}(\Delta t^2)$ when $\nu^w$ and $\nu^v$ are sufficiently smooth). Changes in the moments of the sending domain are 
\begin{align}   
   \Delta \mzero^v &=  - \Mzero^v\;, \label{eq:mzerosendchange}\\
   \Delta \mone^v &= \dela^v \mzero^v(t) - \big(\Mzero^v + \Mone^v \big)\;, \label{eq:monesendchange}\\
   \Delta \mtwo^v &= 2 \dela^v \mone^v(t) + (\dela^v)^2 \mzero^v(t) - \big(\Mzero^v + 2\Mone^v + \Mtwo^v\big)\;, \label{eq:mtwosendchange}
\end{align}
while changes in the receiving domain are 
\begin{align}
   \Delta^v \mzero^w  &= \Mzero^v\;, \label{eq:mzeroreceivechange}\\
   \Delta^v \mone^w &= \gamma^{v \to w} \Mone^v\;, \label{eq:monereceivechange}\\
   \Delta^v \mtwo^w &= (\gamma^{v \to w})^2 \Mtwo^v\;, \label{eq:mtworeceivechange}
\end{align}
where $\Delta^v$ denotes the change due to influx from domain $v$. The easiest change to consider is in $\mzero^v$. The sending domain loses mass in \eqref{eq:mzerosendchange} that gets sent to the receiving domain \eqref{eq:mzeroreceivechange}, thus conserving the total mass. One may compare such an approach with finite volume schemes \cite{LeVeque2002} where the numerical flux leaving a finite volume cell is identical to the flux that is received by an adjacent cell. Naturally, equations \eqref{eq:monesendchange} and \eqref{eq:mtwosendchange} represent changes in the first and second moments of the sending domain due to loss of mass and advection. The changes to the first and second moments of the receiving domain \eqref{eq:monereceivechange} and \eqref{eq:mtworeceivechange} account for the mass influx and any variations in advection speed represented by $\gamma^{v \to w}$. 
In theory, the formulas can never lead to negative moments. However, in reality, minor numerical errors in integration may lead to negative moments on the order of the error in the lookup table or the integration scheme. Thus, in the numerical computation of \eqref{eq:mzerosendchange}--\eqref{eq:mtworeceivechange} we implement a cap on the change of the moments to ensure non-negativity.

Having derived the changes in the moments in the sending and receiving domains of a given network edge, we now collect the total change of moments in a given domain $v \in \mV$, which has both an inflow edge and an outflow edge. From \eqref{eq:mzeroreceivechange}--\eqref{eq:mtworeceivechange} we obtain the $\Delta^u$-terms, caused by the domain $u$ which sends information into domain $v$; and the $\Delta$-terms \eqref{eq:mzerosendchange}--\eqref{eq:mtwosendchange} are caused by the outflux from the domain $v$ itself (which are unaffected by where that flux actually goes). Hence, the resulting update in domain $v$ is:
\begin{equation}
\label{eq:apmethodupdate}
    m^v_i(t+\Delta t) = m^v_i(t) + \Delta^u m^v_i + \Delta m^v_i, \quad \text{for moments~} i \in \{0,1,2\}\;.
\end{equation}

\subsection{Discussion of moment methods in contrast to other numerical methods}
\label{subsec:other_methods}
Conceptually, the network PDE problems described in \S\ref{subsec:governing_equations} are not complicated to numerically approximate via standard finite volume methods \cite{LeVeque2002}, or similar approaches like DG \cite{CockburnKarniadakisShu2000}. Such methods are conservative by construction and their framework is naturally compatible with fluxes across domains. For instance, for traffic flow networks \cite{GaravelloPiccoli2006}, the basic Godunov method generates what in transportation engineering is called ``cell transmission models'' \cite{Daganzo1994}. There are two fundamental complications with such approaches. First, if diffusion is present in \ref{eqn:domain_PDE}, the problem is more stiff than plain advection, thus (semi-)implicit methods may be required to avoid the need to tiny time steps. Second, even for pure advection, fixed-grid methods tend to produce spurious deformations of the numerical solution. For first order methods, that comes (to leading order) in the form of numerical diffusion that results in a spurious smearing of profiles. Higher order methods, when chosen linear, tend to produce spurious overshoots that can yield spurious negativity of the numerical solution. Carefully designed limiters \cite{LeVeque2002} can remedy this last problem, but at the expense of simplicity of the method. Moreover, one obtains nonlinear numerical methods even though the original problem is linear.

For these reasons, for the specific application considered here, fixed-grid methods were found in \cite{lewkiewicz_temperature_2022} to be far from ideal for the transport of very peaked Gaussians caused by the mechanism in \S\ref{subsec:mechanism}. Instead, a moving mesh finite volume method was developed, for which the computational mesh on each domain $v$ moves with the respective advection speed $\nu^v(t)$. In- and out-fluxes through the domain boundaries are carefully accounted for, and the decay and diffusion terms of \ref{eqn:domain_PDE} are then treated on that moving mesh exactly and implicitly, respectively. In the absence of diffusion, the  specialized numerical scheme is exact within each domain; the only approximation occurs for the fluxes between domains. A comparison in \cite{lewkiewicz_temperature_2022} revealed that in some scenarios, the moving mesh method required 100 times fewer grid cells than a Godunov method, despite both methods being first order. Like schemes with limiters, the moving-mesh method is also nonlinear (in the mesh-shift variable). However, as used in \S\ref{sec:results_slf}, it is possible to employ the method to generate one-year-forward mapping operators that are linear by accepting a small amount of interpolation error.

The moving mesh method can (up to the chosen mesh resolution) capture arbitrary solutions on the domain, which renders it a suitable method for the general model study conducted in \cite{lewkiewicz_temperature_2022}. However, for those (far from uncommon) situations where the solution is close to a Gaussian, the method still requires 50--100 cells per domain. This is where the moment methods presented above have their key role: they can represent and track Gaussians exactly; and the results in \S\ref{sec:results_slf} reveal that even for profiles not extremely close to a Gaussian, critical quantities of interest may nevertheless be well-reproduced. The moment methods are nonlinear approximations of a linear problem, a property shared with fixed-grid methods with limiters and moving mesh schemes.

A structural distinction between the two moment methods constructed above is analogous to two ways to discretize space-time PDE: the ODE-based method in \S\ref{subsec:ode_scheme} is the analog of semi-discrete ``method of lines'' approaches that first approximate the space variable only, and then apply suitable time-stepping schemes to the resulting ODE system; in contrast, the asymptotic preserving method in \S\ref{subsec:ap_scheme} is analogous to fully discrete methodologies, as those developed in \cite{LeVeque2002}. In the important regime of strongly peaked solutions, the ODE-based method incurs rapidly varying right hand sides and thus requires tiny time steps. In contrast, the method in \S\ref{subsec:ap_scheme} is asymptotic preserving: because the spatial integrals for very narrow Gaussians are easy to know exactly, no significant time-step restriction is incurred. In fact, while not implemented here, theoretically the method could even represent and track Dirac delta solutions. This is another key advantage of moment methods: mesh-based methods (fixed or moving) have a hard time capturing peaks that are narrower than the mesh resolution.

Many of the features of the moment methods here are also shared with entropy closures \cite{BrunnerHolloway2001}, as used for instance in radiation transport simulations \cite{ChidyagwaiFrankSchneiderSeibold2018}: nonlinear, potential to represent strongly peaked solutions (radiation beams), positivity. And also the fact that---except for simple cases---the closure, respectively reconstruction, cannot be written analytically and thus requires numerical quadrature or lookup tables.

\vspace{1.5em}
\section{Numerical Results: Showcasing Success/Failure on Key Test Problems}
\label{sec:results_test_problems}
This section highlights the performance (and challenges) of the developed moment methods on a carefully designed simple test problem: the two-domain advection-only problem. The density is initialized as a Gaussian centered at $a = \frac{1}{2}$ in domain~1, and constant zero in domain~2. The total mass is always $1$. The initial value of $\sigma$ is varied, and the initial second moment $\mathring{m}^1_2$ is calculated to reproduce that Gaussian width. The advection speed is taken to be constant $\nu = 1$ across both domains, and the inflow for domain~1 is the outflow of domain~2, and vice-versa. Thus the true solution to the PDE has the following moments:
\begin{equation*}
(m^1_0,m^1_1,m^1_2,m^2_0,m^2_1,m^2_2)(t) =
\begin{cases}
  (1,\frac{1}{2},\mathring{m}_2,0,0,0)\quad & t=2n   \text{ for }n=\{0,1,2\dots\} \\
  (0,0,0,1,\frac{1}{2},\mathring{m}_2)      & t=2n+1 \text{ for }n=\{0,1,2\dots\}\;.
\end{cases}
\end{equation*}
While this is the behavior of the moments of the true solution of the advection PDE, it is not the true solution of the moment ODE. Hence, this simple test case can be used to demonstrate under which circumstances the solution of the moment ODE is close to the moments of the true solution of the PDE.

\subsection{Long-time behavior of moment methods}
Here, we numerically investigate the long time behavior of the ODE-based moment method applied to the above two-domain advection-only problem. The ODE is discretized in time with RK4, and flux limiting as discussed in \S\ref{subsec:ode_scheme}. In particular, we consider the case when the time step is sufficiently small, i.e., $\Delta t < \sigma$, and numerical time stepping errors are negligibly small.

\begin{figure}
    \includegraphics[width=.7\linewidth]{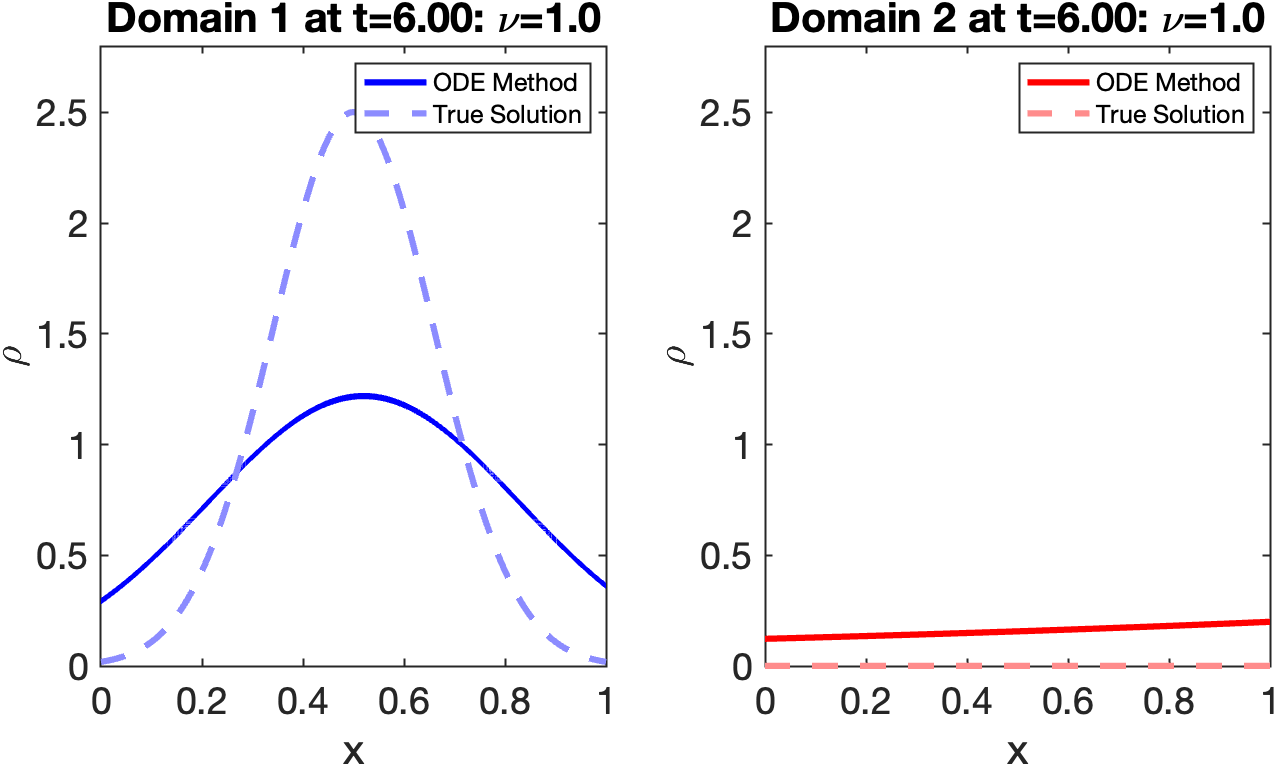}
    \caption{Two-domain advection test, used for 5.1 and 5.2. This particular example visualized the spurious widening of the Gaussians, happening when the solution fails to be well-localized ``within a single domain''. When the solution is well-localized, the moment method would be visually exact. The true solution to the PDE is shown dashed in red and blue.}
    \label{fig:gaussian_spread}
\end{figure}

Errors may occur when the exact solution is not realizable by a single Gaussian at every time, or if the true solution is not well approximated by the Gaussian recovered from the moments. If significant mass is both entering and leaving the same domain at the same time (which may occur should the Gaussian be sufficiently wide), on that domain the recovered Gaussian will exhibit smearing, and the outflux will be decreased, causing a spurious widening of the Gaussian in the next domain as well, see Fig.~\ref{fig:gaussian_spread}. This effect will compound over time, and the flow will eventually be driven to a constant. Because Gaussians are not compactly supported, this effect will in principle occur no matter the width of the initial Gaussian. However, if the Gaussian is well localized, the effect will be completely negligible.

\begin{figure}
    \includegraphics[width=.50\linewidth]{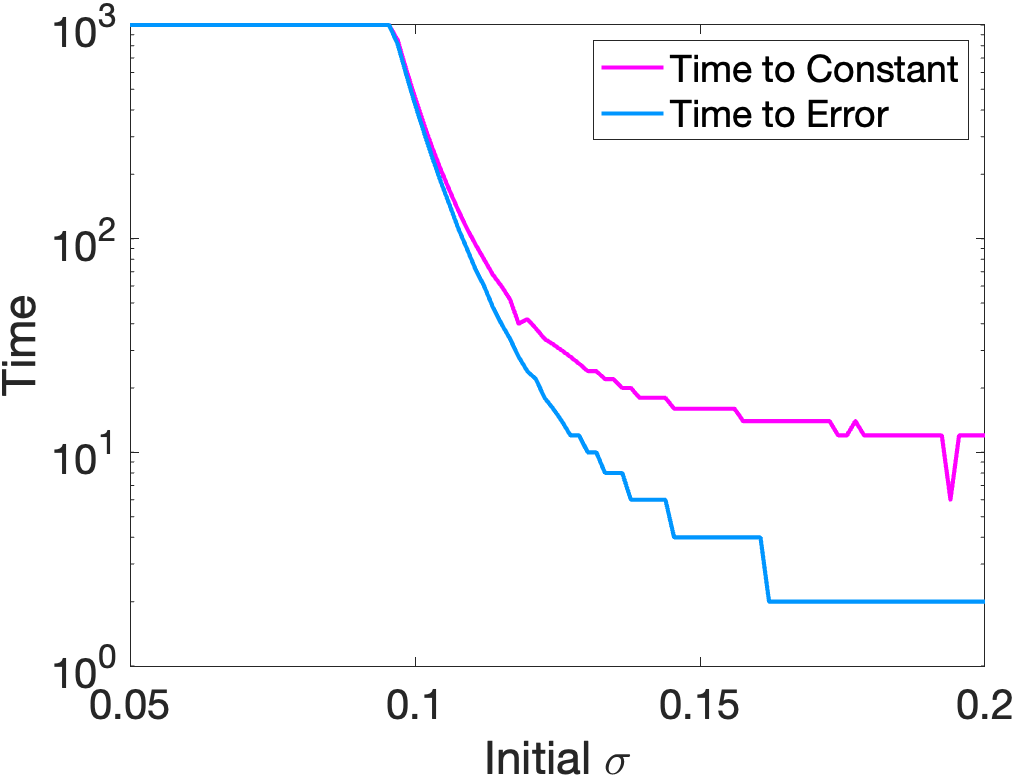}

    \caption{Long time behavior of ODE-based moment method, RK4 time step $\Delta t = 0.005$. Error is calculated at $t=2n$ for $n\in\mathbb{N}$. The time for the density to become constant is shown in magenta. The time for the density to a relative error of $0.01$ is shown in blue. Final time is capped at $t=1000$. Wide Gaussians are driven to constant due to the smearing effects seen in Fig.~\ref{fig:gaussian_spread}. Sufficiently narrow but well-resolved Gaussians do not spread significantly even over large time scales.}
    
    \label{fig:MOL_error}
\end{figure}

Figure~\ref{fig:MOL_error} examines this effect based on initial $\sigma$-value. We do not examine the behavior of extremely narrow Gaussians in this plot, as the ODE-based method is not designed to perform well for such problems. For wide Gaussians, we see the flow driven to constant quickly. This occurs when the true solution cannot be well represented by a single Gaussian. In this case, the moment approach is fundamentally ill-suited to the problem (we note that this effect will also be seen for the asymptotic preserving approach studied in \S\ref{sec:results_need_for_AP}). However, for well-localized Gaussians, a spurious increase in variance is not seen, i.e., the solution is not driven to constant in at least 500 cycles ($t = 1000$). We further emphasize that for these well-localized Gaussians, we maintain a relative error in the moments of less than $0.01$. This demonstrates that although the reconstruction mapping may introduce small errors, the errors do not compound, and accuracy is maintained over long time scales. 

On the note of the spurious widening of the Gaussians, it is intriguing to compare this with standard finite volume methods. In the Godunov method, numerical diffusion will produce a similar effect. However, it would occur in the opposite fashion: numerical diffusion would be particularly dominant for strongly peaked Gaussians, while the moment method does particularly well. Moreover, achieving no significant deformation of a finite volume solution over 500 cycles would require an extremely fine grid resolution.

\subsection{Need for an asymptotic preserving moment method for peaked solutions}
\label{sec:results_need_for_AP}
As we have seen, the ODE-based moment method incurs spurious widening for the simple two-domain advection in the case where $\sigma$ is large. On the other side of the spectrum, we now consider the case when $\sigma$ is small, i.e., we have peaked Gaussians as described in \S\ref{subsec:mechanism}. In this case, the asymptotic preserving (AP) moment method from \S\ref{subsec:ap_scheme} is natural. A narrow Gaussian is captured via three moments, and being well localized, the spurious spread as in Fig.~\ref{fig:gaussian_spread} does not occur. 
However, as discussed in \S\ref{subsec:ap_scheme}, when the time scale $\tau$ is much smaller than the time step $\Delta t$, we will observe large numerical errors in our calculation of the moments, in particular with an explicit time stepping scheme. The AP scheme is designed to avoid this, and to produce accurate results even when the time step is large. For the specific problem at hand, the AP method from \S\ref{subsec:ap_scheme} reads as follows. For domain~1 we have
\begin{align*}
\mzero^1(t+\Delta t)&= \mzero^1(t) + \Delta \mzero^1 + \Delta^2 \mzero^1\;, \\
\mone^1(t+\Delta t)&= \mone^1(t) + \Delta \mone^1  + \Delta^2 \mone^1\;, \\
\mtwo^1(t+\Delta t)&= \mtwo^1(t) + \Delta \mtwo^1 + \Delta^2 \mtwo^2\;,
\end{align*} 
and the update in domain~2 is
\begin{align*}
   \mzero^2(t+\Delta t) &= \mzero^2(t) + \Delta \mzero^2 + \Delta^1 \mzero^2\;, \\
   \mone^2(t+\Delta t)&= \mone^2(t) + \Delta \mone^2 + \Delta^1\mone^2\;, \\
   \mtwo^2(t+\Delta t)&= \mtwo^2(t) + \Delta \mtwo^2 + \Delta^1 \mtwo^2\;.
\end{align*}

To test the accuracy of the AP method versus the ODE-based method for very narrow Gaussians, we examine the relative error in $\sigma$ after one cycle, or t=2. At t=2, the Gaussian has returned to its initial position, and we can easily extract the error in $\sigma$. Both methods are mass-conservative, and errors in the expected value are found to be negligible (at least 5 order of magnitude smaller than the error in $\sigma$). Thus, tracking the error in the width gives a clear picture of when the methods succeed or when they fail.

\begin{figure}
    \includegraphics[width=.32\linewidth]{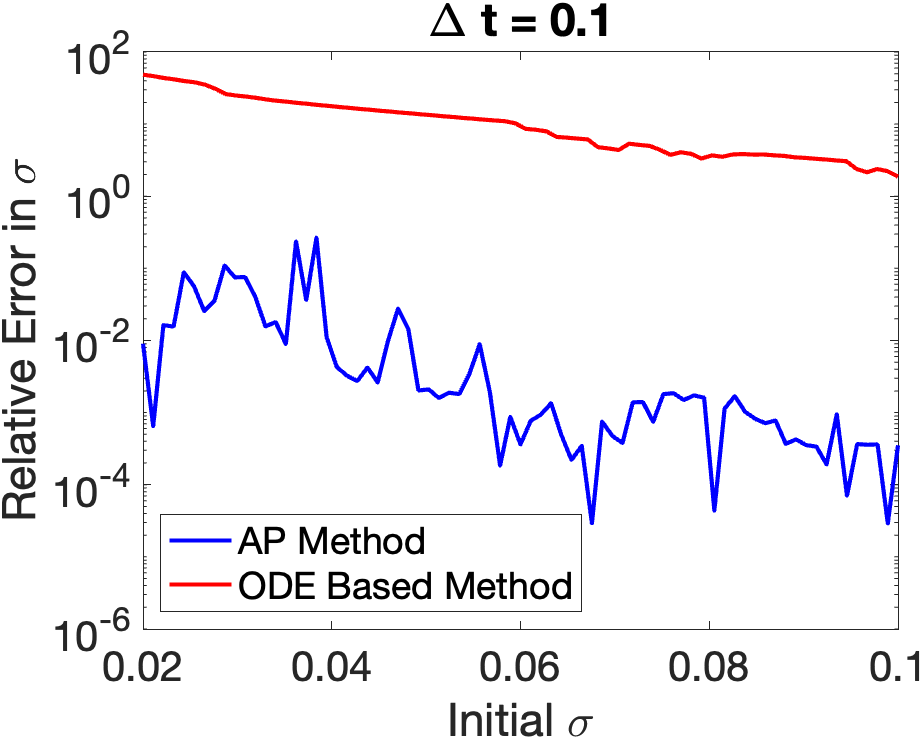}
    \includegraphics[width=.32\linewidth]{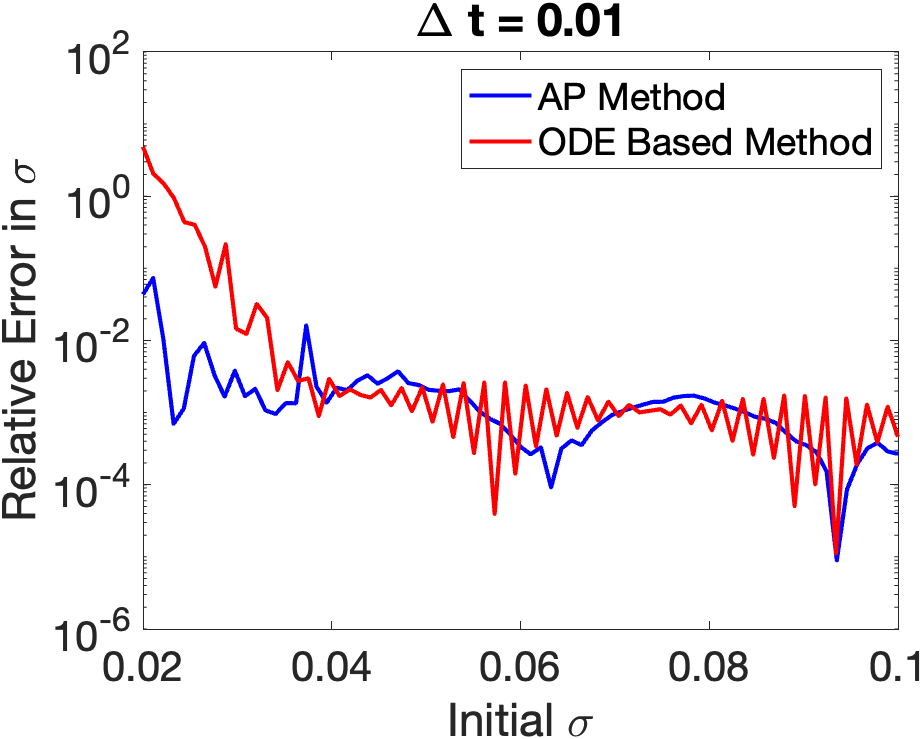}
    \includegraphics[width=.32\linewidth]{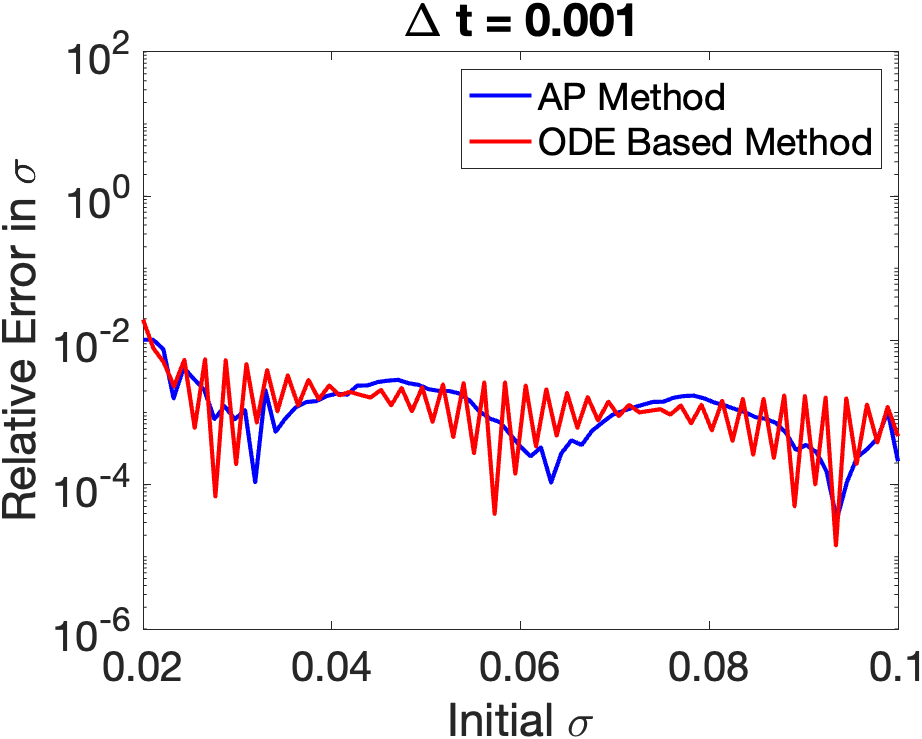}
    \caption{The relative error in $\sigma$ at t=2 (one cycle) is shown for three choices of time step size. For the largest time step size, only the AP method, shown in blue, is accurate. For a moderate time step size, we see the AP method is necessary for narrow Gaussians, but the ODE-based methods, shown in red, perform similarly for wider Gaussians. Finally, for a very small time step size the methods perform comparably. Apparent noise in the errors is dominated by the error in the reconstruction mapping. }
    \label{fig:APvsMOL}
\end{figure}

Figure~\ref{fig:APvsMOL} demonstrates that for narrow Gaussians, the AP method accurately recovers the correct $\sigma$ value even when the time step is large. Indeed, the error in the AP method applied to this simple problem appears nearly independent of the size of the time step. As expected, we need small time steps to resolve narrow Gaussians using the ODE-based method. Though the ODE-based approach is accurate for small time steps, for large time steps compared to the width of the Gaussian, the error is large.

\begin{remark}
Though errors in Fig.~\ref{fig:APvsMOL} appear noisy, this is because they are dominated by the error in the reconstruction. We note that for these tests, $\sigma$ is only expected to be accurate to four significant digits. As $\Delta t \to 0$, we see the error approach the error in the lookup table for both methods. With sufficiently small time steps, the error in both methods reduces to solely the reconstruction error.
\end{remark}



\vspace{1.5em}
\section{Application Results: Establishment Maps for Invasive Forest Pests}
\label{sec:results_slf}
We now apply the ODE-based moment method (\S\ref{subsec:ode_scheme}) to the specific network PDE model that describes the life cycle of spotted lanternfly (SLF). As described in \S\ref{subsec:application_slf} and visualized in Fig.~\ref{fig:network_slf}, the SLF model deviates from the generic model framework (from \S\ref{subsec:governing_equations}) in one aspect: the influx into egg domains, defined in \eqref{eq:slf_model_egg_laying}, is determined via a weighted integral over the motile density. We compute this integral via a simple quadrature rule, employing precisely the same reconstructed density $\rho^v(a)$ as used for the other terms in \eqref{eq:moment_system_ode}.

\begin{figure}
    \includegraphics[width=1\linewidth]{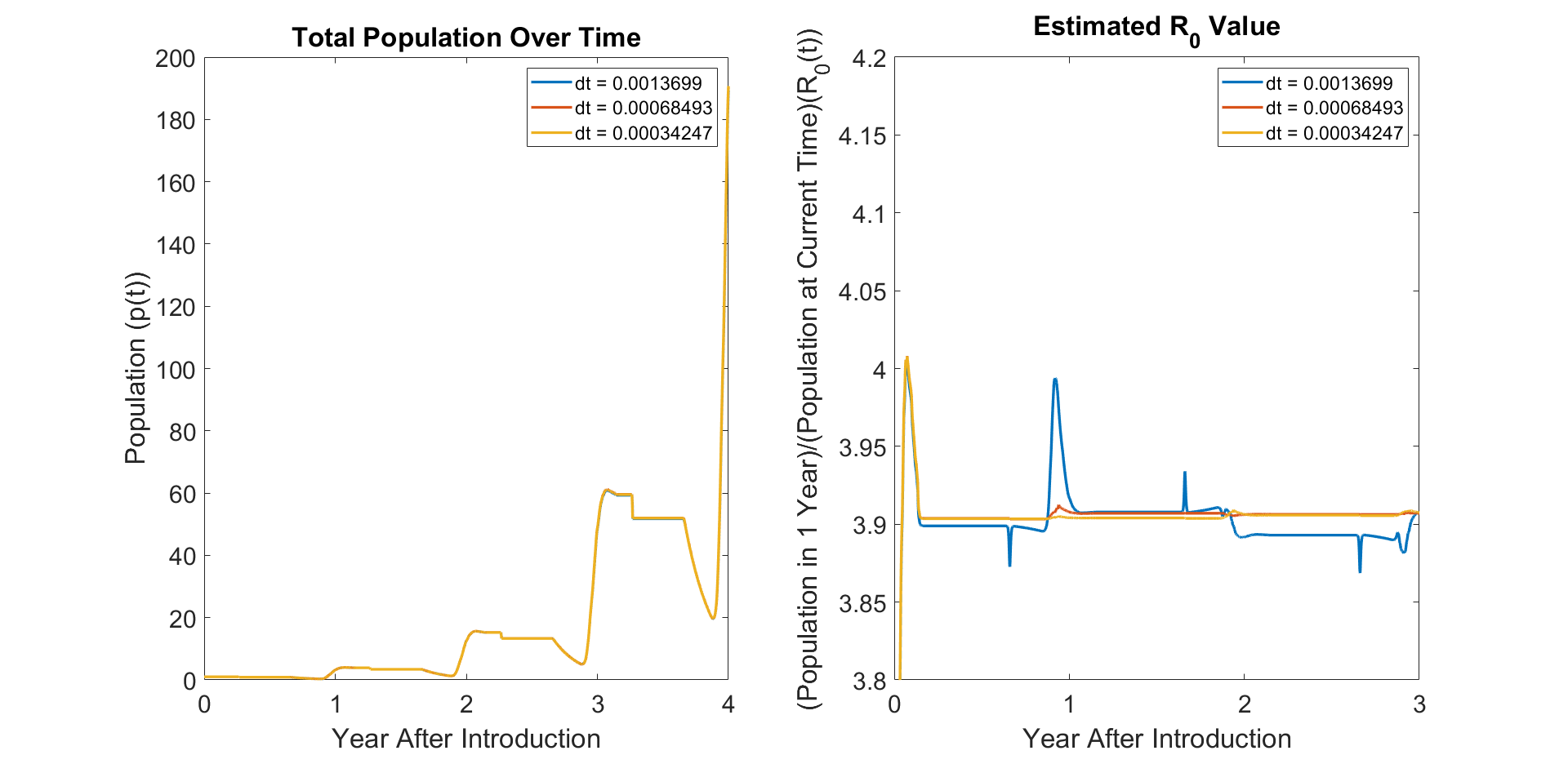}
    \caption{Numerical study of time step needed. Shown is the total population $p(t)$ (left) and derived $R_0(t)$ (right), computed using the ODE-based moment method (\S\ref{subsec:ode_scheme}), integrated via RK4.}
    \label{fig:totalpopulationrnaught}
\end{figure}

Previous studies of the full SLF life cycle model \cite{lewkiewicz_temperature_2022} revealed that wherever diapause synchronizes the population dynamics with the annual cycle, the population distribution in each domain tends to be close to a single Gaussian. In these situations, the moment method is expected to be able to accurately capture population dynamics. In contrast, when such a synchronization does not occur, the population distributions may not be captured well by a single Gaussian and thus the moment method cannot be expected to be accurate. The effects of diapause depend on the annual temperature cycle to which the population is exposed.

\begin{figure}
    \includegraphics[width=1\linewidth]{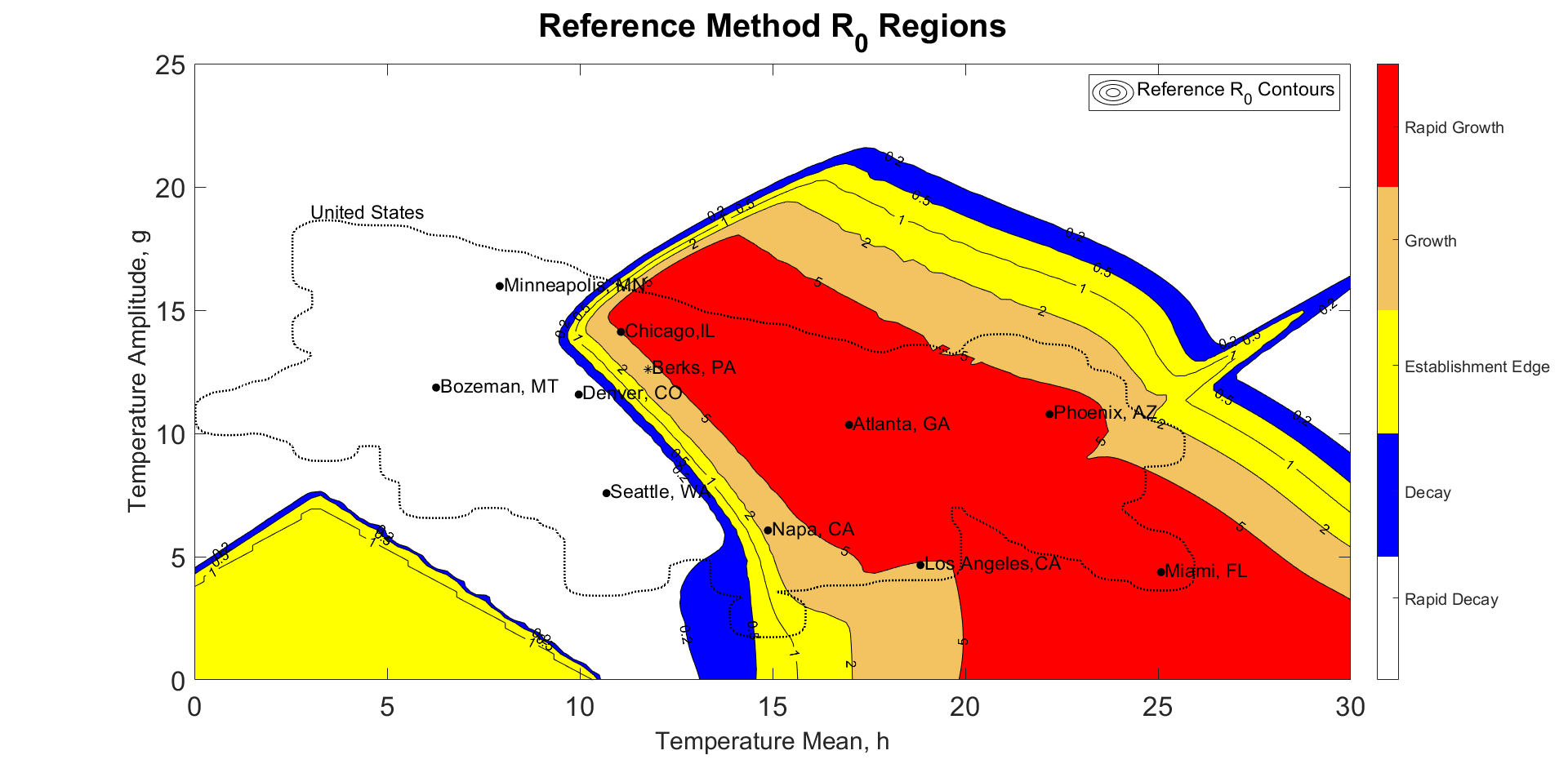}
    \includegraphics[width=1\linewidth]{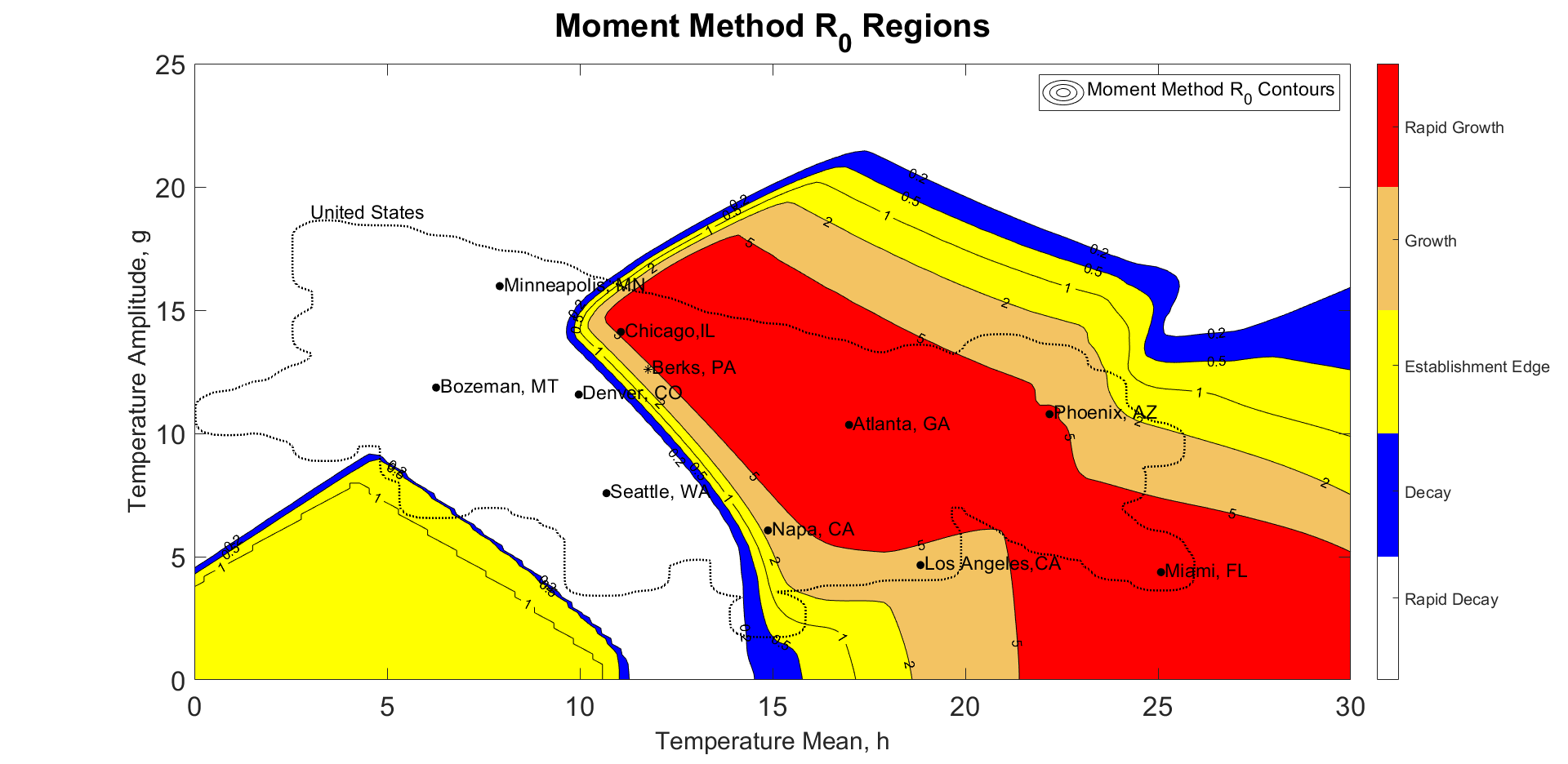}
    \caption{Population development regions, based on reproductive number $R_0$, corresponding to different levels of growth. Overlaid is the border of the US in the parameter space and some large cities and the initial invasion location, Berks County.}
    \label{fig:parametersweep}
\end{figure}

In order to understand how well the moment method performs in various temperature regimes we study the reproductive number $R_0$ of the population, here defined as the one-year growth factor. The $R_0$-value is an important quantity in population dynamics (and other applications like disease modeling) like those of invasive forest pests. For populations synchronized with the annual cycle, $R_0$ measures how many offspring a single individual produces on average, and it indicates whether a population is able to persist in a given location. Values $R_0 > 1$ correspond to a growing population while values $R_0 < 1$ indicate shrinking populations. 

Using the reference moving mesh method (see \S\ref{subsec:other_methods}) on Eqn.~\eqref{eqn:domain_PDE}, one can formulate the linear one-year solution operator and then use its dominant eigenvalue, the maximum achievable one-year growth, as a proxy for the $R_0$-value. All reference solutions and reference $R_0$-values herein are produced via the moving mesh method.

\begin{figure}
    \includegraphics[width=1\linewidth]{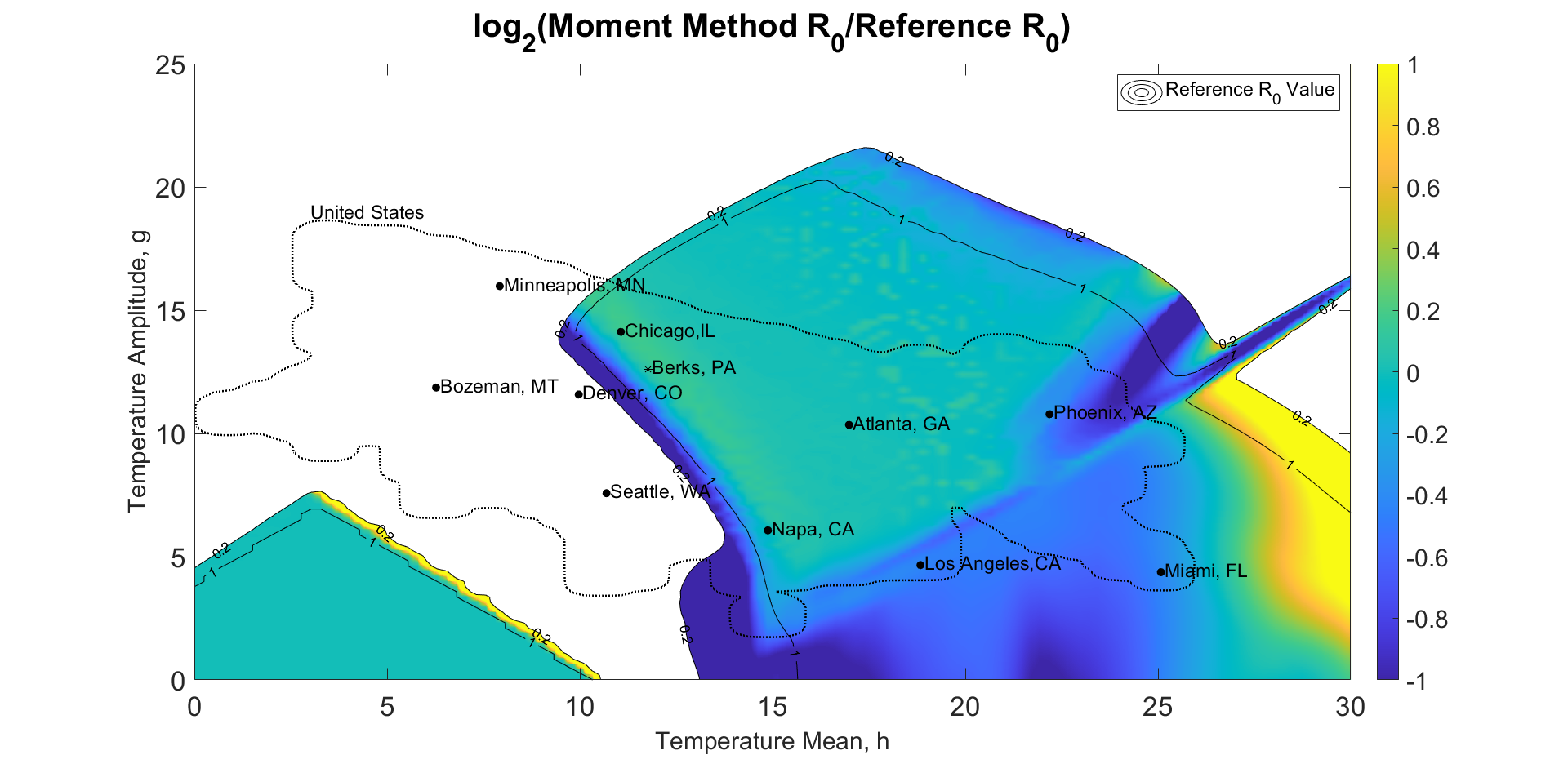}
    \caption{Comparison of the $R_0$-value obtained from the simulations and the Reference $R_0$-value for $R_0$-values above $0.2$.}
    \label{fig:R0error}
\end{figure}

Because the moment method is nonlinear (see \S\ref{subsec:other_methods}), instead of formulating a linear eigenvalue problem, the $R_0$-value is instead obtained via forward simulations over multiple years, analogous to a power iteration to find a dominant eigenvalue. We consider the total population across all SLF life stages
\begin{equation*}
p(t) = \sum_{v=1}^4 \int_{0}^{1} \rho^v(a,t) \ud{a}\;,
\end{equation*}
and calculate the $R_0$-value at a given time by the formula $R_0(t) = \frac{p(t+1)}{p(t)}$. For all simulations we initialize populations to have all individuals in the diapause domain with $E=0.04975$ and $V=0.0475$; however, the choice of initial condition does not affect the resulting $R_0$-value significantly. 
The simulation results in Fig.~\ref{fig:totalpopulationrnaught} demonstrate that after an initial transient phase during the first year, the $R_0$ estimator quickly convergences to a constant function. Even for small times (less than a year), $R_0(t)$ is close to constant in time. Large deviations occur when there are large changes in the population e.g., during egg laying. The $R_0$-value is estimated as a temporal average of $R_0(t)$.

Figure~\ref{fig:totalpopulationrnaught} shows on the left an example of how a population that experiences the effects of diapause grows over time. Due to the diapause process we see a stagnation of the population growth (early in the year). After the eggs hatch and before the motiles are mature, we see a small decrease in population due to cold death. Once the motile population begins to lay eggs again we return to the stage of rapid growth. Here we see an initial population of 1 egg, which over the course of four years grows to a population of 200 individuals. We emphasize that this simplified model is linear, and does not include integer effects, thus this factor of 200 may be taken as a multiplier for any initial population. Many areas of Pennsylvania have seen rapid infestation, with even a small population taking hold quickly. 

\begin{figure}
    \includegraphics[width=.90\linewidth]{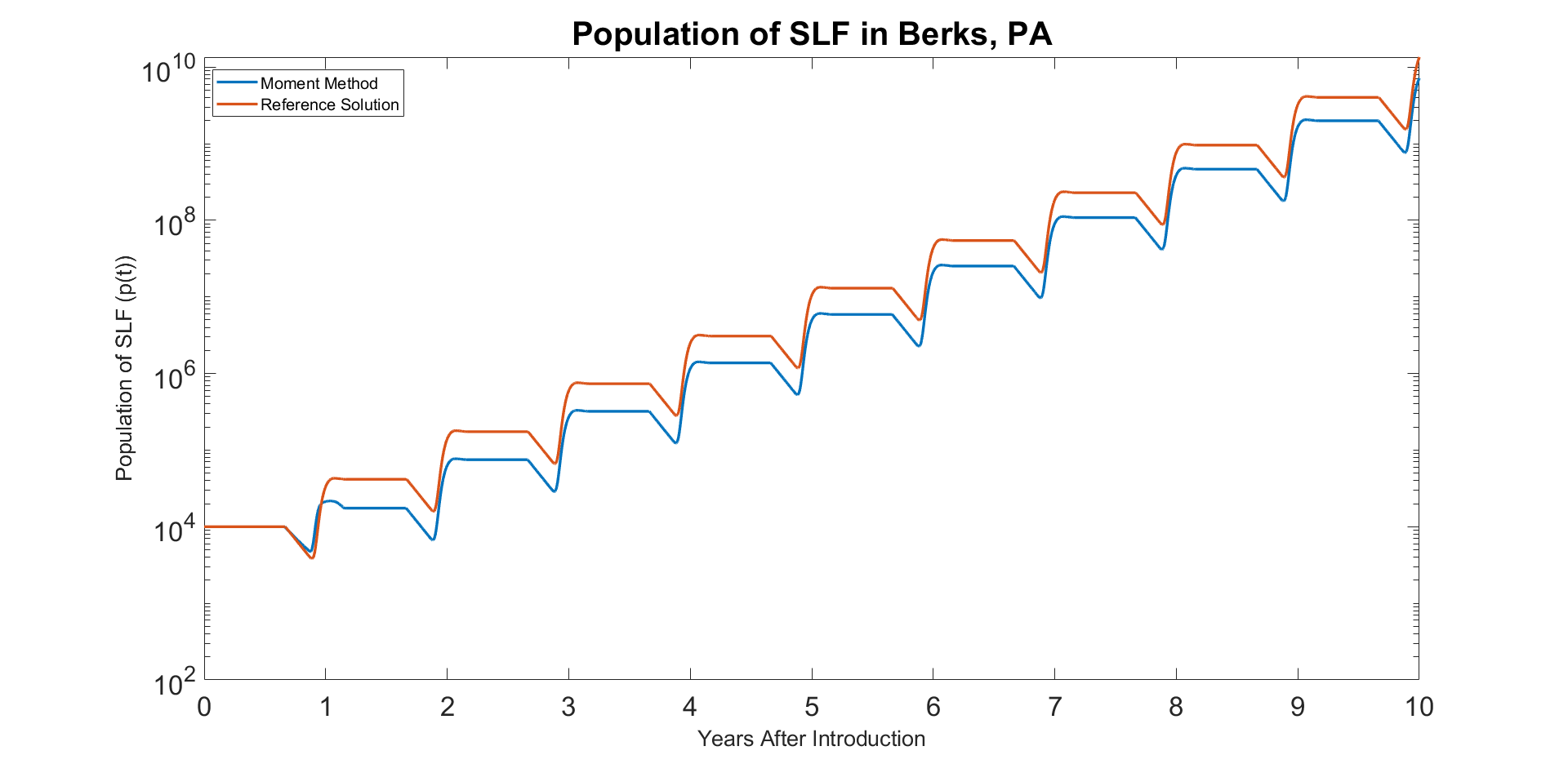}
    \includegraphics[width=.90\linewidth]{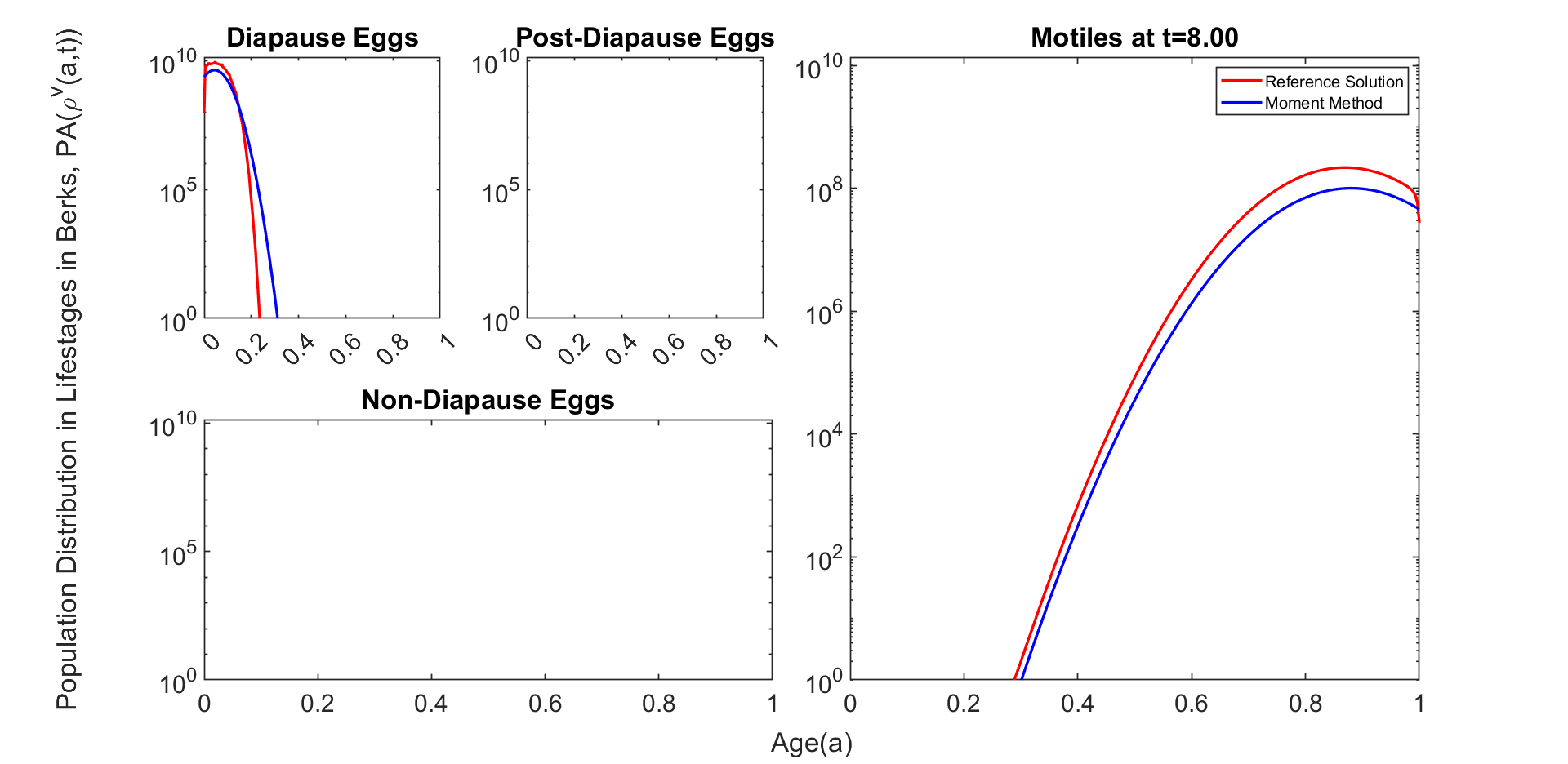}
    \caption{SLF populations for Berks, PA, the initial location of the infestation in the US. The top shows the total population over time $p(t)$. The bottom shows the densities for each life stage $\rho^v(a,8)$, i.e., 8 years after infestation.}
    \label{fig:slf_results_Berks}
\end{figure}

\begin{figure}
    \includegraphics[width=.90\linewidth]{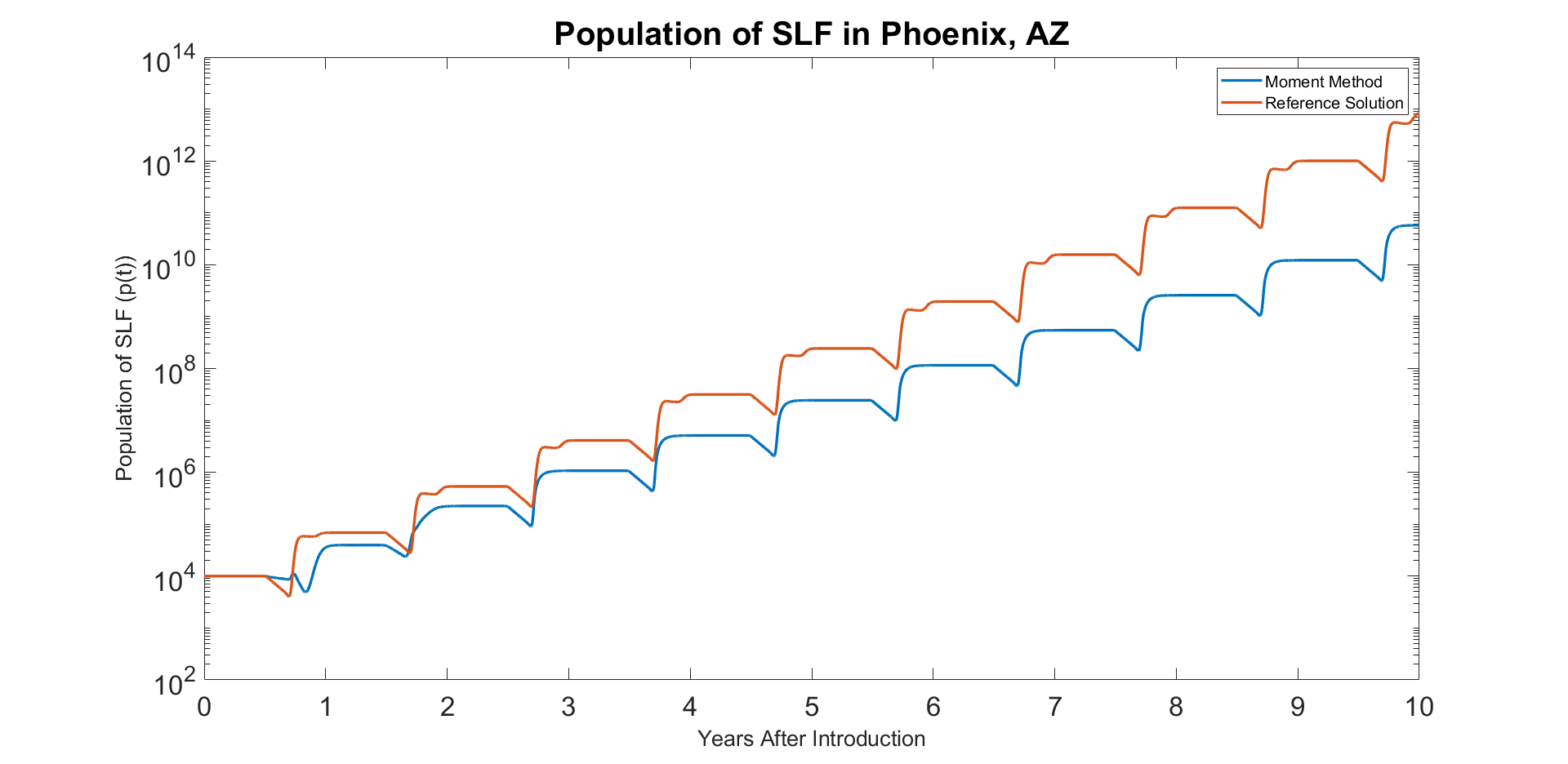}
    \includegraphics[width=.90\linewidth]{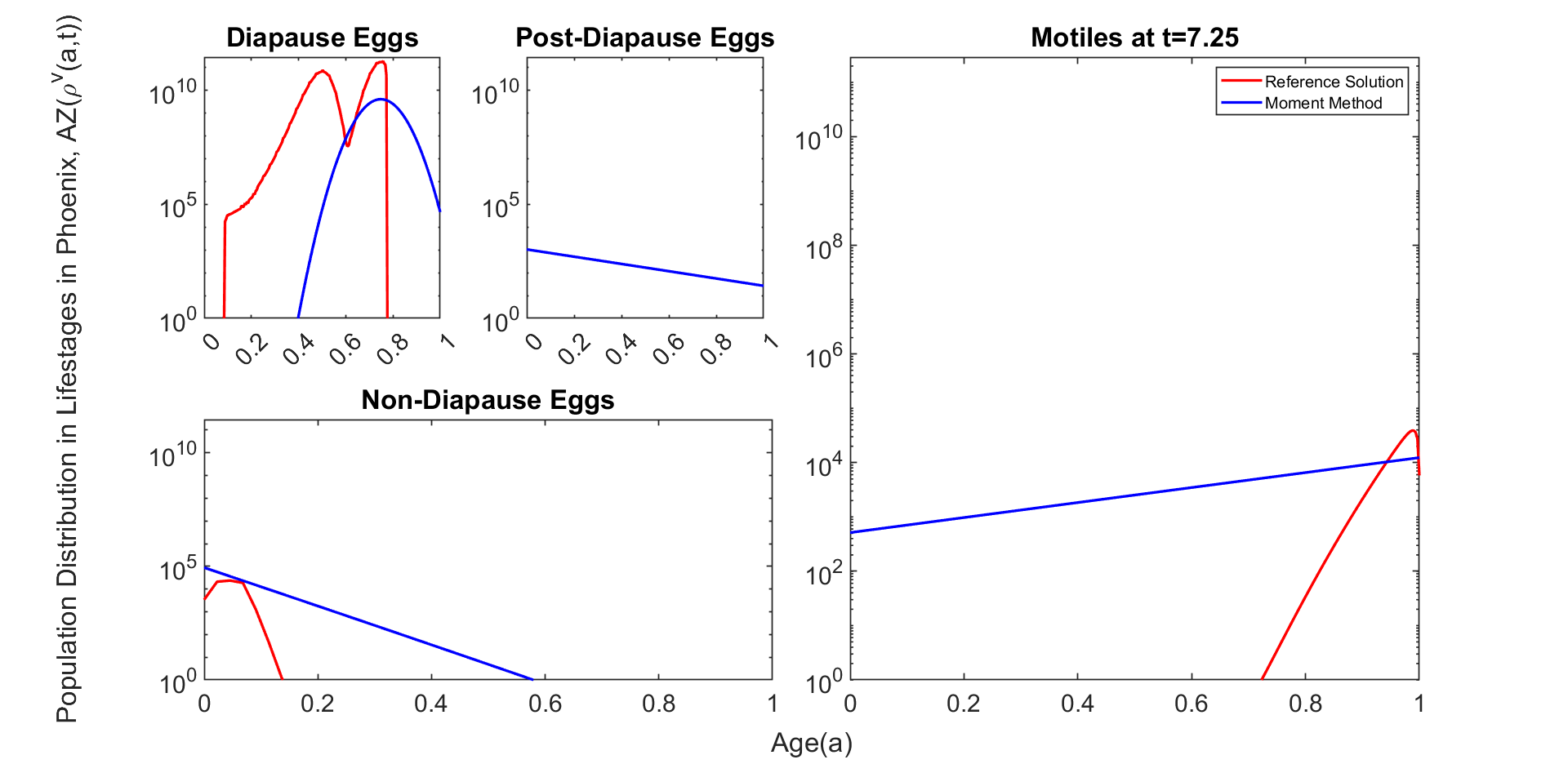}
    \caption{SLF populations for Phoenix, AZ. The top shows the total population over time $p(t)$. The bottom shows the densities for each life stage $\rho^v(a,7.25)$, i.e., 7.25 years after infestation. The bimodal densities are clearly not captured by the moment method, but the trend of rapid growth is seen for both reference and moment solutions.}
    \label{fig:slf_results_Phoenix}
\end{figure}

\begin{figure}
    \includegraphics[width=.90\linewidth]{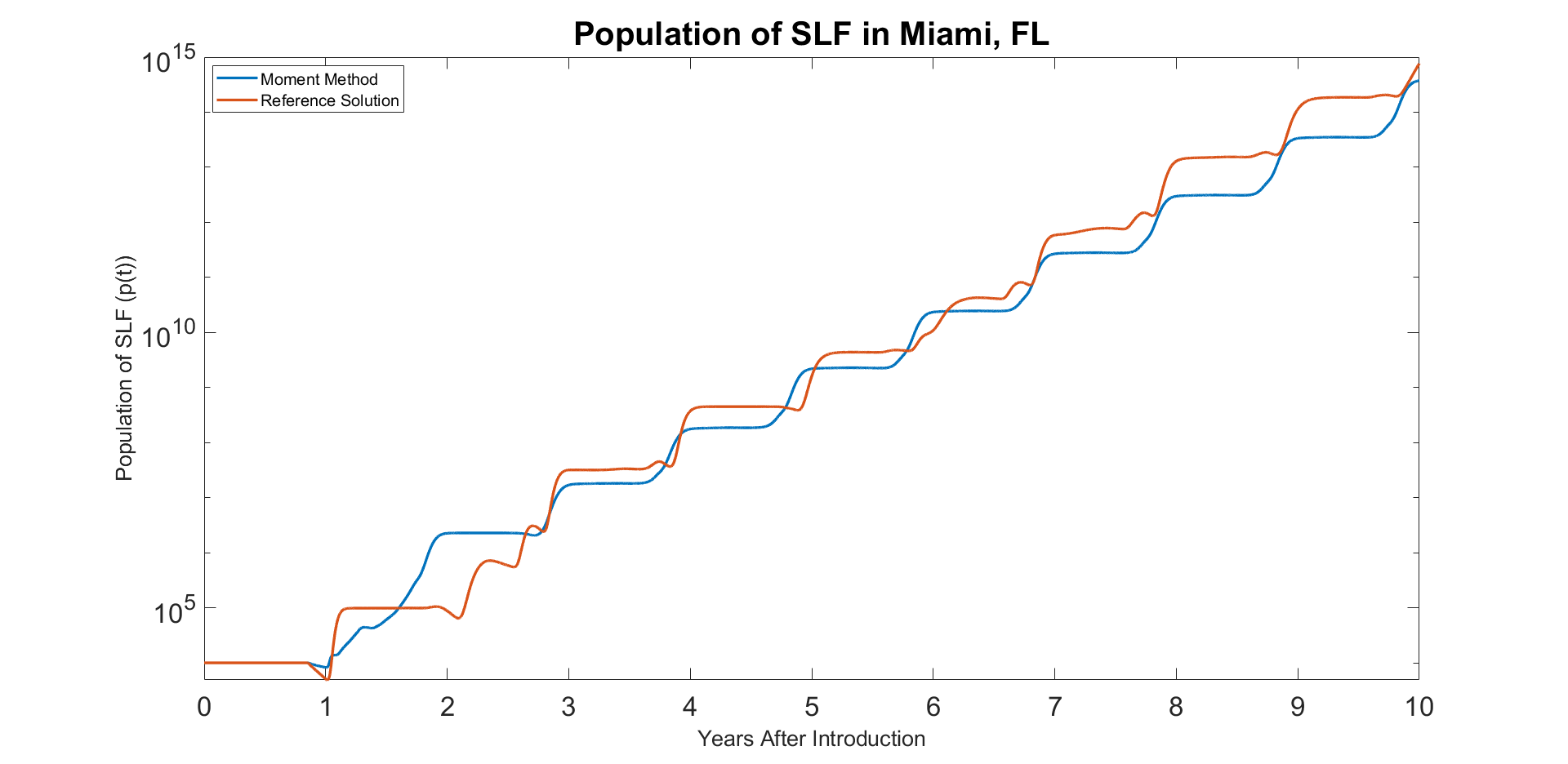}
    \includegraphics[width=.90\linewidth]{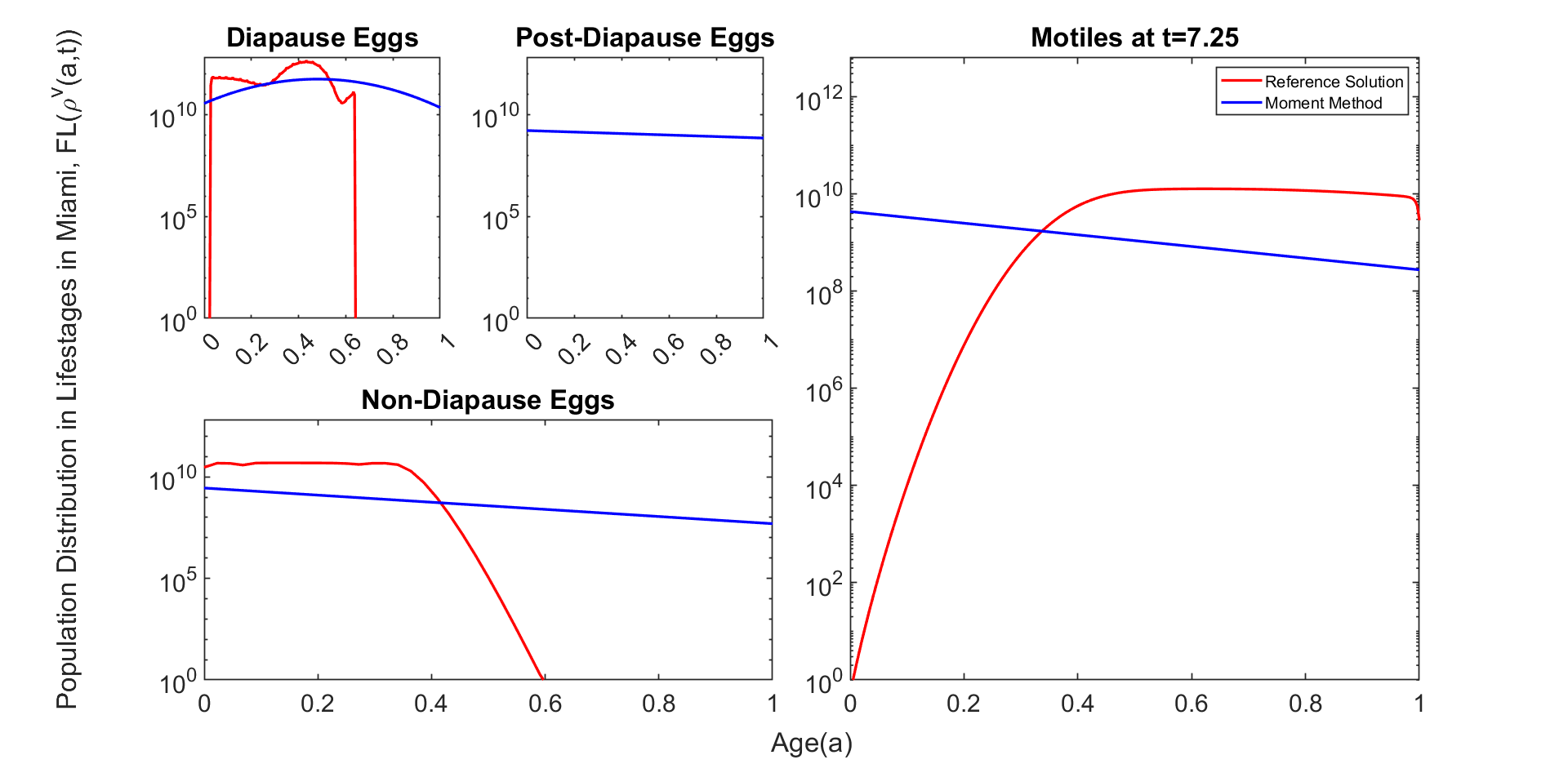}
    \caption{SLF populations for Miami, FL. The top shows the total population over time $p(t)$. The bottom shows the densities for each life stage $\rho^v(a,7.25)$, i.e., 7.25 years after infestation. Even though the densities are not well reproduced, the resulting $R_0$ estimate is very good.}
    \label{fig:slf_results_Miami}
\end{figure}

To assess how the computed $R_0$-values from the moment method compare to the reference $R_0$-values, we consider temperature profiles that are sinusoidal with mean $h$ (average annual temperature) and amplitude  $g$ (temperature variation), i.e.,
\begin{equation*}
T(t) = g\sin(2\pi t)+h\;,
\end{equation*}
where $t = 0$ corresponds to April 1 \cite{lewkiewicz_temperature_2022}. This simplified temperature profile is chosen to capture the long-term establishment potential in a given location under average conditions devoid of noise. In the $(h,g)$-parameter space a full study is conducted, simulating each temperature profile for 4 years, with an initial population of 1000 eggs, all in the diapause life stage with the same $(E,V)$-values given above.

To identify the risk of SLF establishment, we define $R_0$-value ranges. An $R_0<0.2$ is considered ``rapid decay''; $R_0 \in [0.2,0.5)$ is denoted ``decay''; $R_0\in [0.5,2)$ is the ``establishment edge''; $R_0\in [2,5)$ is categorized as ``growth''; and $R_0\ge 5$ is categorized as ``rapid growth''.
In Fig.~\ref{fig:parametersweep} we see our parameter space broken into the 5 different regions. The $R_0$-value contours are plotted for both the reference method and the moment method. For easy orientation, the plot is overlaid with the border of the United States and some large cities. When comparing the $R_0$ values of the moment method against the reference $R_0$, we see a good general agreement between the two. In some temperature regimes, diapause does not synchronize development, thus the moment method is unlikely to be able to accurately capture the population density. However, we do not see a significant degradation in the agreement of the methods even for those temperature regimes. A notable region is the triangle in the bottom left of both plots. This corresponds to a temperature regime where diapause eggs neither develop nor die. The contour of this regime is well captured with the moment method. In contrast, the high-temperature mid-amplitude regime is not well captured by the moment method. In particular, the distinctive spike extending the establishment edge is lacking. In these temperature regimes there can be multiple SLF generations per year. As discussed in \S\ref{sec:results_test_problems}, the moment method cannot accurately capture multiple peaks in one domain. Instead, the resulting population curve will be flattened, and outflow will be decreased, leading to inaccurate results that propagate to other domains. 

Figure~\ref{fig:R0error} shows the phase-space map of the ratio between the estimated $R_0$-value and the reference $R_0$. Regions where the reference $R_0 < 0.2$ are masked out, because in them a large relative $R_0$ error does not change the prediction of rapid population decay. For high mean temperatures the moment method does not do well at estimating the $R_0$-value. This is expected due to the lack of diapause synchronization. To highlight and understand the specific shortcomings of the moment method, we study three selected reference points by running numerical simulations over 10 years: (A) Berks County, PA, the initial area the invasive species was found in the US and a region where we expect the moment method to work well; (B) Phoenix, AZ, where a comparatively large error is observed; and (C) Miami, FL, where the synchronization effects of diapause do not manifest.

The results for Berks County, PA, shown in Fig.~\ref{fig:slf_results_Berks}, demonstrate that indeed the moment method works very well: $p(t)$ agrees nicely with the reference solution, aside from a slight undershoot in population count after the first year; however, this shift does not affect the resulting $R_0$-value. We also see that the population distribution 8 years after the introduction is well captured by a single Gaussian in each domain. 

For Phoenix, AZ, shown in Fig.~\ref{fig:slf_results_Phoenix}, we see that the reference population has some features that are lost in the moment method. The reference total population curve has two periods per year of rapid growth, in contrast to the moment method's single phase of rapid growth. The distributions for each life stage indicate the reason. The reference solution has a bimodal density in the diapause domain, which cannot be reproduced by the moment method. The bimodal distribution indicates that there is a bit of a synchronization effect of diapause but it is not strong enough to render the population distribution unimodal. We also see that all other domains have a nonzero population, indicating the overlap of multiple developmental pathways at the same time. While such types of solutions are not what the moment method is designed to reproduce, the resulting $R_0$ estimate is relatively accurate; both the reference solution and the moment method classify the region as one of rapid growth.

For Miami, FL, shown in Fig.~\ref{fig:slf_results_Miami}, similar aspects as with Phoenix are present: the high temperatures and low temperature variability allow for the superposition of multiple developmental pathways. Consequently, the densities in each domain are not reproduced well by a single Gaussian. Again, the trajectory of the total population \emph{is} again reproduced well by the moment method.


\vspace{1.5em}
\section{Conclusion and Outlook}
\label{sec:conclusions}
This study demonstrates that for certain applications of network PDE problems, suitably designed moment methods can significantly reduce the dimensionality of the problem, and at the same time capture the relevant structure and behavior of the solution satisfactorily in a wide range of problem parameters. To some extent, the success of the moment framework does rely on the fact that the network PDE problem has a tendency to produce solutions that are close to Gaussians, and that such solutions can be well captured via only three moments on each domain. However, the numerical results of the SLF application in Phoenix (Fig.~\ref{fig:slf_results_Phoenix}) and in Miami (Fig.~\ref{fig:slf_results_Miami}) demonstrate that even in cases where the true solution is not well described by a single Gaussian (and thus the moment method cannot reproduce the true solution well), relevant quantities of interest (here: the reproductive number $R_0$) can still be captured within an acceptable margin of error (see Fig.~\ref{fig:R0error}).

The methods presented herein rely on a proper pre-computation of the mappings between moments and reconstructions (\S\ref{subsec:moment_reconstruction}). Once these are established, the simplest resulting approach is an ODE-based moment method (\S\ref{subsec:ode_scheme}). It transforms the multi-domain PDE problem into system of ODEs, whose dimension is as small as three times the number of domains. For the SLF application, that means that a 12-dimensional ODE system results. If the reconstructed Gaussians are not too peaked, the resulting ODE system can be solved via standard ODE solvers, like Matlab's \texttt{ode45.m}, or preferably: methods that ensure positivity (see \S\ref{subsec:ode_scheme}). Compared to numerical methods that directly discretize the original multi-domain PDEs \cite{LeVeque1992}, the ODE-based moment method represents a significant structural simplification, achieving high-order in time is easier, and the computational complexity is much lower.

For problems that exhibit strongly peaked solutions, as it is frequently the case in the absence of diffusion, the ODE method is not ideal, see \S\ref{sec:results_need_for_AP}. For such situations, the asymptotic preserving (AP) moment method presented in \S\ref{subsec:ap_scheme} captures the correct solution behavior without requiring the time step to resolve the fast time scales induced by moving narrow Gaussians. In fact, while not implemented here, there is no limitation in the moment method from allowing the tracking of solutions that possess Dirac delta peaks.

An important property of the moment models is that they are nonlinear, even though the original network PDE problem is linear. That fact incurs some challenges, for instance the loss of a superposition principle (see the $R_0$ construction in \S\ref{sec:results_slf}). That being said, the approximation of a linear PDE problem via a nonlinear numerical scheme is not at all uncommon; for instance, high-order TVD schemes for advection require nonlinearity \cite{LeVeque1992}. The reconstruction employed here is also not dissimilar to entropy closures. These also approximate the linear radiation transport equation via nonlinear models \cite{BrunnerHolloway2001}.

For the specific application of predicting establishment potential of spotted lanternfly in the United States, the results in \S\ref{sec:results_slf} reveal that moment methods may indeed represent an intriguing modeling and computational asset. First, they reduce a multi-domain PDE model to an ODE system, which is both a conceptual and a computational simplification. And second, the results are largely intriguing. In temperature regimes where development is fully synchronized with diapause, the moment method accurately captures both the approximate growth rate of the population and the population density on each domain. In contrast, where diapause synchronization is not dominant, the reference density profiles are frequently not captures well. However, even in those situations, the approximate growth rate of a population is largely reproduced well nevertheless.

Based on the methods and results established herein, the following future work aspects arise naturally. First, while the moment methods are demonstrated to perform well in many important scenarios, important numerical analysis questions remain, most prominently: do the methods guarantee the dynamic conservation of realizability?

Another important direction is how to further improve the methods' efficiency. Rather than pre-computing a lookup table of reconstructions, one could---specific to the application---pre-compute and store directly the quantities needed in the moment method, like values at outflow boundaries or the egg-laying kernel in the SLF application. Such efficiency improvement can allow truly rapid simulations of the model, as needed when embedded into a larger framework like optimization or optimal control. In fact, for optimization, even just the structural simplification of having an ODE instead of a multi-domain PDE could be beneficial. Another important pursuit are methods that guarantee the methods' robustness, such as exploring suitable time-stepping schemes for the moment methods that guarantee conservation and positivity, such as \cite{BurchardDeleersnijderMeister2003}.

The AP method has shown to address the challenge of strongly peaked solutions satisfactorily for pure advection; hence, its extension to problems involving diffusion and decay is a natural goal. Another important question on the generalization of the methods is: which of the concepts established here could be carried over to tracking more than three moments per domain? In that situation, the reconstruction problem becomes significantly more general, and entropy-based \cite{tagliani_numerical_2001} or spline-based techniques \cite{john_techniques_2007, lebaz_reconstruction_2016} could be explored.

On the modeling of population dynamics, a natural direction is to extend the model to incorporate spatial spread \cite{JUNG2017, LADIN2023}. The model would then become a multi-domain problem in $(x,y,a,t)$, where $(x,y)$ is the geo-spatial location, while $(a,t)$ are the same variables as here. For such a model framework, the dimensional reduction to $(x,y,t)$, that is offered by moment methods, can be critical, both towards the applicability of existing software and towards computational efficiency.


\vspace{1.5em}
\section*{Acknowledgments}
The authors would like to thank Stephanie Lewkiewicz for assistance with the simulation software for the calibrated spotted lanternfly model.
This material is based upon work supported by the National Science Foundation under Grant No.~DMS--2309728 (Seibold) and Grant No.~DMS--2202888 (Kean). Any opinions, findings, and conclusions or recommendations expressed in this material are those of the authors and do not necessarily reflect the views of the National Science Foundation.
Research was sponsored by the DEVCOM Analysis Center and was accomplished under Cooperative Agreement Number W911NF-22-2-0001. The views and conclusions contained in this document are those of the authors and should not be interpreted as representing the official policies, either expressed or implied, of the Army Research Office or the U.S. Government. The U.S. Government is authorized to reproduce and distribute reprints for Government purposes notwithstanding any copyright notation herein.

\vspace{1.5em}
\bibliographystyle{abbrv}
\bibliography{focused_references,references_seibold}

\vspace{2.5em}
\end{document}